\RequirePackage{fix-cm}
\documentclass[smallextended]{svjour3}
\smartqed  
\usepackage{graphics}
\usepackage{epsfig}
\usepackage{amssymb}
\usepackage{amsfonts}
\usepackage{amsmath}
\usepackage{graphicx}
\begin{document}
\title{On the piecewise pseudo almost periodic solution of nondensely impulsive integro-differential systems with infinite delay}
\author{
Adn\`ene Arbi and Farouk Ch\'erif}
\institute{
Adn\`ene Arbi\\
$adnen.arbi@enseignant.edunet.tn, arbiadnene@yahoo.fr$\\
National School of Advanced Sciences and Technologies of Borj Cedria, University of Carthage, El Khawarizmi Street, Carthage 2078, Tunisia;
Laboratory of Engineering Mathematics (LR01ES13), Tunisia
Polytechnic School, University of Carthage, El Khawarizmi Street, Carthage 2078, Tunisia.\\
Farouk Ch\'erif\\
$faroukcheriff@yahoo.fr.$\\
University of Sousse, Department of Computer
science, ISSATS, Laboratory of Math Physics; Specials Functions and
Applications LR11ES35, Ecole Sup\'{e}rieure des Sciences et de
Technologie, Sousse 4002, Tunisia.
   }
\maketitle
\begin{abstract} In the theory of neutral differential equations with pulse influence
(neutral impulsive differential equations), there are many unsolved
problems related to certain results in the theory of integral and
integro-differential equations. In this article, we present an
result for the existence of the piecewise pseudo almost periodic
solution of a class of neutral impulsive nonlinear
integro-differential systems with infinite delay. The method used
involves result on the theory of integrated semigroup as well as the
Sadovskii's fixed point theorem.   
\keywords{Impulsive integro-differential equations, Pseudo almost periodic
solution, Sadovskii's fixed point theorem, Infinitesimal generator
of an analytic semigroup.}
\end{abstract}        



\section{Introduction}
Because they can describe various of real processes and phenomena,
the delayed impulsive differential equations is an important theory
in mathematical modeling of phenomena, we can cited for example in:
artificial neural networks, biologies, computer sciences and physics
(\cite{abada}, \cite{acta1}, \cite{aacta},
\cite{neuro}, \cite{liu}, \cite{sak}). We refer the reader to the
monographs of Bainov and Simeonov \cite{bainov}, Lakshmikantham et
al. \cite{bbainov},  and Samoilenko and Perestyuk \cite{samoil}. In
\cite{shan}, Li Shan Liu et al. obtain a unique solution of the
first order impulsive integrodifferential equation of mixed type in
a Banach space, an explicit iterative approximation of the solution
and estimate the approximation sequence. Using the
Krasnoselski-Schaefer type fixed point theorem in \cite{chang}, the
authors prove the existence of solutions for impulsive partial
neutral functional differential equations with infinite delay in a
Banach space. In \cite{zabut}, Stamov et al. employ the contraction
mapping principle to proved some sufficient conditions for the
existence of almost periodic solutions of impulsive
integro-differential neural networks.

On the other hand, many researchers are interested to periodic
concept, in particular studied the existence and stability of
periodic solutions (\cite{henrikez}, \cite{liz}, \cite{zang}) and
the references cited therein. However, upon considering long–term
dynamical behaviors it is possible for the various components of the
model to be periodic with rationally independent periods, and
therefore it is more reasonable to consider the various parameters
of models to be changing almost-periodically rather than
periodically. Thus, it is more reasonable to consider the almost
periodic behavior of solutions. The investigation of almost periodic
solutions is established to be more accordant with reality. Although
it has widespread applications in real life, the generalization to
the notion of almost periodicity is not as developed as that of
periodic solutions. To the best of authors' knowledge, there are a
few recent published papers considering the notion of almost
periodicity of differential equations with impulses (\cite{sstamov},
\cite{zabut}).

Motivated by the works mentioned above and noting that there is no
papers in the literature for pseudo almost periodic solution of
nondensely impulsive neutral functional integrodifferential
equations with infinite delay, the main purpose of this paper is to
establish the existence and of the piecewise pseudo almost periodic
(PPAP) solution of an impulsive integro-differential equations in a
general Banach space modeled in the form:
\begin{eqnarray}\label{elsevier1}
&&\frac{d}{dt}\left[x(t)-g\left(t,x_{t},\int_{0}^{t}h(t,s,x_{s})ds\right)\right]= A\left[x(t)-g\left(t,x_{t},\int_{0}^{t}h(t,s,x_{s})ds\right)\right]\nonumber\\
&+& f\left(t,x_{t},\int_{0}^{t}k(t,s,x_{s})ds\right),\,\ \text{if}\,\ t\in J=[0,b]\backslash t_{k},\,\ k=1,...,m, \\
&& \label{elsevier11}\Delta
x(t_{k})=x\left(t_{k}^{+}\right)-x\left(t_{k}^{-}\right)=I_{k}\left(x\left(t_{k}^{-}\right)\right),\,\ \text{if}\,\ t=t_{k},\\
&&\label{elsevier111} x_{0}=\phi\in \mathrm{BM}_{h},
\end{eqnarray}
where\\

(i) $A:D(A)\subset X\longrightarrow X$ is a non-densely defined
($\overline{D(A)}\neq X$) closed linear operator and satisfies the
Hille-Yosida condition (see Definition \ref{hy}), (where $X$ is a
Banach space with the norm $\|\cdot\|_{X}$).

(ii) The nonlinear functions $f,g:J\times \mathrm{BM}_{h} \times
X\longrightarrow X$, $h,k:\Delta\times
\mathrm{BM}_{h}\longrightarrow X$, $\mathrm{BM}_{h}$ is the abstract
phase space which will be defined later,
$$
\Delta=\{(t,s): \,\ 0\leq s\leq t\leq b\},
$$
$x(t_{k}^{+})$ and $x(t_{k}^{-})$ represent the right and left
limits of of the function $x$ at $t_{k}$, respectively.

(iii) $I_{k}\in C(X,X)$ $(k=1,...,m,\,\ \text{where}\,\
m\in\mathbb{Z}_{+})$ are bounded functions.

(iv) For any $x:(-\infty,b]\longrightarrow X$ and $t\in J$, the
histories $x_{t}:(-\infty,b]\longrightarrow X$ defined by
$x_{t}(s)=x(t+s)$, $s\leq0$ which belong to $\mathrm{BM}_{h}$.\\

It is well known that the solution $x(\cdot)$ of the problem
(\ref{elsevier1})-(\ref{elsevier111}) is a piecewise continuous
function with points of discontinuity at the moments $t=t_{k}$,
$k\in\mathbb{Z}$, at which it is continuous from the left, i.e. the
following relations are valid:
$$
x(t_{k}^{-})=x(t_{k}),
$$
and
$$
x\left(t_{k}^{+}\right)-x\left(t_{k}^{-}\right)=I_{k}\left(x\left(t_{k}^{-}\right)\right),
\,\ \forall k\in\mathbb{Z}.
$$

For examples of the linear operators with nondense domain satisfying
the Hille-Yosida condition, we can look in the one-dimensional,
 $E=C\left([0,1],\mathbb{R}\right)$ and defined the
operator $A:D(A)\longrightarrow E$ by $Ay'=y,$ where $
D(A)=\left\{y\in C^{1}\left([0,1],\mathbb{R}\right):\right.$\\
$\left.y(0)=0\right\}.$ Then $ \overline{D(A)}=\left\{y\in E:
y(0)=0\right\}\neq E. $ In the n-dimensional, let
$\Omega\subset\mathbb{R}^{n}$ be a bounded open set with regular
boundary $\Gamma$ and define $E=C(\overline{\Omega},\mathbb{R})$ and
the operator $A:D(A)\longrightarrow E$ defined by $Ay=\Delta y$,
where $ D(A)=\left\{y\in E: y=0\,\ \text{on} \,\ \Gamma; \Delta y\in
E\right\}. $ Here $\Delta$ is the Laplacian in the sense of
distributions on $\Omega$. In this case we have $
\overline{D(A)}=\left\{y\in E: y=0 \,\ \text{on} \,\
\Gamma\right\}\neq E.$ See the paper found by Da Prato and
Sinestrari \cite{prato}, Abada et al. \cite{abada} for more examples
and details concerning
the non-densely defined operators.\\

In the case where $g(\cdot)$ is a null function, $f:J\times
\mathrm{BM}_{h} \longrightarrow X$, and $A$ is a nondensely closed
defined linear operator generating a $C_{0}$-semigroup of bounded
linear operators, the problem (\ref{elsevier1})-(\ref{elsevier111})
has been investigated on compact intervals in \cite{abada}. To the
best of our knowledge, no paper in the literature has investigated
the existence, uniqueness and globally exponential stability of PPAP
solution for system (\ref{elsevier1})-(\ref{elsevier111}). Hence,
our goal in this paper is to study the dynamics of model
(\ref{elsevier1})-(\ref{elsevier111}). By applying Sadovskii's
theorem, we give some sufficient conditions ensuring the existence
of PPAP solution of system (\ref{elsevier1})-(\ref{elsevier111}),
which are new and complement the previously known results.\\

The remainder of this paper is organized as follows: In Section 2,
we will introduce an abstract phase space $\mathrm{BM}_{h}$ and the
concept of integrated semigroups. In Section 3, the concept of
pseudo almost periodicity is presented, we derived also, some
preliminaries and hypotheses which will be used in the paper. In
section 4, we present a new idea of research to prove some criteria
for ensuring the existence of the PPAP solution. At last, an
illustrative example is given. It should be mentioned that the main
results include Theorem \ref{th1}. The results of the present paper
extend to a nondensely defined operator some ones considered in the
previous literature.

Our results are based on the properties of the analytic semigroup
and ideas and techniques in Kellermann
\cite{kellerman}, Pazy \cite{pazy}, Yosida \cite{yosida}.\\

Let $\mathbb{R}$ and $\mathbb{Z}$ be the sets of real and integer
number. Define the space $\mathbf{T}$ by $\mathbf{T}=\{\{t_{k}\}:
t_{k}\in(-\infty,+\infty), t_{k}< t_{k+1},
k=\pm1,\pm2,...,\lim\limits_{k\longrightarrow\pm\infty}t_{k}=\pm\infty
\}$, we denote the set of unbounded and strictly increasing
sequences. Let $A:D(A)\longrightarrow X$ be the infinitesimal
generator of a compact analytic semigroup of bounded linear
operators $\{S(t)\}_{t\geq0}$ on a Banach space $X$ with the norm
$\|\cdot\|$, and let $0\in\rho(A)$, then it is possible to define
the fractional power $(-A)^{-\alpha}$, for $0\leq \alpha<1$, as
closed linear invertible operator with domain $D((-A)^{-\alpha})$
dense in $X$. The closedness of $D((-A)^{-\alpha})$ implies that
$D((-A)^{-\alpha})$ endowed with the graph norm
$\|x\|_{X}=\|x\|+\|(-A)^{-\alpha}x\|$ is a Banach space. Since
$(-A)^{-\alpha}$ is invertible, its graph norm $\|x\|_{X}$ is
equivalent to the norm $|x|=\|(-A)^{-\alpha}x\|$. Thus
$D((-A)^{-\alpha})$ equipped with the norm $|\cdot|$ is a Banach
space which we denote by $X_{\alpha}$.\\

To obtain the existence of the pseudo almost
periodic solutions of (\ref{elsevier1})-(\ref{elsevier111}), we make the following assumptions:\\\\
($H_{1}$) For every $t\in (-\infty,0]$ the function
$t\longrightarrow \phi(t)$ is pseudo almost periodic, and for $t\in
J$, the functions $t\longrightarrow
f\left(t,x(\cdot),y(\cdot)\right)$, $t\longrightarrow
g\left(t,x(\cdot),y(\cdot)\right)$ are pseudo almost periodic for
$x(\cdot),y(\cdot)$ pseudo almost periodic and there exist constants
$G_{i}>0$, ($i=1,2$) and $H>0$ such
that\\

(i) $\|g(t,\psi,x)-g(t,\chi,y)\|_{X}\leq
G_{1}\|\psi-\chi\|_{\mathrm{BM}_{h}}+G_{2}\|x-y\|_{X}, \,\ t\in J,
\,\
\psi,\chi\in\mathrm{BM}_{h}, \,\ x,y\in X$;\\

(ii) $\|h(t,s,\psi)-h(t,s,\chi)\|_{X}\leq
H\|\psi-\chi\|_{\mathrm{BM}_{h}}, t,s\in J, \,\
\psi,\chi\in\mathrm{BM}_{h}$.\\

(iii) $\|(-A)^{\beta}f(t,x,\psi)-(-A)^{\beta}f(t,y,\chi)\|_{X}\leq
F_{1}\|x-y\|_{\mathrm{BM}_{h}}, t,s\in J, \,\
\psi,\chi\in\mathrm{BM}_{h}$,\,\ \text{where}\,\ $\beta\in(0,1)$ and
$F_{1}>0$.
\\\\
($H_{2}$) (i) Let $A$ be an infinitesimal operator of an analytic
semigroup $S(t)$ and that $0\in\rho(T)$, such that
$$
\left\|(-A)^{1-\alpha}S(t-s)\right\|\leq
\frac{C_{1-\alpha}}{(t-s)^{1-\alpha}}\exp(-\lambda(t-s)), \,\ t\in
J, \,\ \alpha\in(0,1), \,\ \lambda>0.
$$
%

Denote by $BC(\mathbb{R},X_{\alpha})$ the space of
$X_{\alpha}$-valued bounded continuous functions.\\
(ii) The semigroup $\{S'(t)\}_{t\geq0}$ is pseudo almost periodic,
compact on $(\overline{D(A)},\|\cdot\|)$, whenever $t>0$, and there
exists
$M\geq1$ such that $\|S'(t)\|\leq M$, for all $t\in J$.\\\\
($H_{3}$) (i) For every $t\in\mathbb{R}$, the functions
$t\longrightarrow h(t,s,x_{s})$, $t\longrightarrow k(t,s,x_{s})$ are
pseudo almost periodic for $(s,x_{s})\in\mathbb{R}\times X$ and
there exist continuous functions $p,q:J\longrightarrow[0,+\infty)$
such that
$$
\left\|\int_{0}^{t}h(t,s,\psi)ds\right\|_{X}\leq p(t)\|\psi\|_{C},
\,\ \,\ \left\|\int_{0}^{t}k(t,s,\psi)ds\right\|_{X}\leq
q(t)\|\psi\|_{C},
$$
for every $(t,s)\in\Delta$ and $\psi\in C$.\\
(ii) $\sup\limits_{t\in J}p(t)<\bar{p}$, and $\sup\limits_{t\in
J}q(t)<\bar{q}$, where $\bar{p}$ and $\bar{q}$ are positive
constants.\\\\
($H_{4}$) For every positive integer $k$ there exits $\alpha_{k}\in
L^{1}(J,[0,+\infty))$ such that
$$
\sup\limits_{\|\psi\|_{C},\|x\|\leq k}\|f(t,\psi,x)\|_{X}\leq
\alpha_{k}(t), \,\ \forall t\in J,
$$ and
$$
\liminf\limits_{k\longrightarrow+\infty}\frac{1}{k}\int_{0}^{b}\alpha_{k}(s)ds=\mu<\infty.
$$
($H_{5}$) (i) For each $k$, ($k=1,...,m$), $I_{k}\in
C(X,\overline{D(A)})$ is pseudo almost periodic sequence,
$$\|I_{k}(x)\|_{\overline{D(A)}}\leq L_{k}(\|x\|_{X}),$$ and
$$
\sum\limits_{0<t_{k}<t}L_{k}\left(\left\|\left(x(t_{k}\right)\right\|_{X}\right)<\infty,\,\
\liminf\limits_{\sigma\longrightarrow+\infty}\frac{L_{k}(\sigma)}{\sigma}=\lambda_{k}<\infty,
$$
where $L_{k}:[0,+\infty)\longrightarrow[0,+\infty)$ is a continuous
nondecreasing function.\\

(ii) $\exists L_{I}>0, \,\
\|I_{k}(x)-I_{k}(y)\|_{\overline{D(A)}}\leq
L_{I}\|x-y\|_{X}.$\\\\
($H_{6}$) (i) For each $(t,s)\in\Delta$, the function
$k(t,s,\cdot):\mathrm{BM}_{h}\longrightarrow X$ is continuous and
for each $\psi\in \mathrm{BM}_{h}$, the function
$k(\cdot,\cdot,\psi):\Delta\longrightarrow X$ is strongly
measurable.\\

(ii) For each $(\psi,x)\in C\times X$, the function
$f(\cdot,\psi,x):J\longrightarrow X$ is strongly measurable.\\

(iii) The set of sequences $\left\{\tau_{k+j}-\tau_{k}\right\}$,
$k,j\in\mathbb{Z}$, $\left\{\tau_{k}\right\}\in\mathbf{T}$ is
uniformly pseudo almost periodic and
$\inf\limits_{k}\left\{\tau_{k+1}-\tau_{k}\right\}=\varrho>0$.\\

It is well known that $BC(\mathbb{R},X_{\alpha})$ is a Banach space
when it is equipped with the sup norm defined by
$$
\|x\|_{\alpha,\infty}:=\sup\limits_{t\in\mathbb{R}}\|x(t)\|_{\alpha},
\,\ \text{for}\,\ x\in BC(\mathbb{R},X_{\alpha}),$$ where the
$\alpha$-norm defined by
$$
\|u\|_{\alpha}:=\|(-A)^{-\alpha}u\|_{X},\,\ \text{for}\,\
u\in\mathcal{D}((-A)^{-\alpha}).
$$
\section{Abstract phase space $\mathrm{BM}_{h}$ and integrated semigroups}
\subsection{Phase space $\mathrm{BM}_{h}$}
Assume $h:(-\infty,0]\longrightarrow[0,+\infty)$ is a continuous
function with\\ $l=\int_{-\infty}^{0}h(s)ds<\infty,$ besides for any
$a>0$, define
$$
\mathrm{BM}_{a}=\left\{\psi:[-a,0]\longrightarrow X \,\ \text{such
that}\,\ \psi(t)\,\ \text{is bounded and measurable}\right\},
$$ and equip the space $\mathrm{BM}_{a}$ by the norm
$$
\|\psi\|_{[-a,0]}=\sup\limits_{[-a,0]}\|\psi(s)\|_{X}, \,\ \forall
\psi\in \mathrm{BM}_{a}.
$$
Let us define
\begin{multline*}
\mathrm{BM}_{h}=\left\{\psi:(-\infty,0]\longrightarrow X \,\
\text{such that for any}\,\ c>0\,\
\psi|_{[-c,0]}\in\mathrm{BM}_{a}\right.\\
\left.\text{and}
\int_{-\infty}^{0}h(s)\|\psi\|_{[s,0]}ds<+\infty\right\}.\phantom{+++++}
\end{multline*}
If $\mathrm{BM}_{h}$ is endowed with the norm,
$$
\|\psi\|_{\mathrm{BM}_{h}}=\int_{-\infty}^{0}h(s)\|\psi\|_{[s,0]}ds,
\,\ \forall \psi\in \mathrm{BM}_{h},
$$
 then it is clear that
$\left(\mathrm{BM}_{h},\|\cdot\|_{\mathrm{BM}_{h}}\right)$ is a
Banach
space.\\

Now, for all $\sigma>0$ we consider the space
\begin{multline*}
\mathrm{BM}_{h}'=\left\{x:(-\sigma,b]\longrightarrow X, \,\
\text{such that}\,\ x_{k}\in C(J_{k},X)\,\ \text{is a piecewise
continuous function}\right.\\
 \text{with points of discontinuity} \,\
\tilde{t}\in[-\sigma,0] \,\ \text{at which}\,\ x(\tilde{t}^{-}) \,\
\text{and} \,\ x(\tilde{t}^{+})\,\ \text{exist} \\
\left. \text{and} \,\ x(\tilde{t}^{-})=x(\tilde{t}), \,\
k=0,1,...,m, \,\
x_{0}=\phi\in\mathrm{BM}_{h}\right\},\phantom{+++++++++}
\end{multline*}
where $x_{k}$ is the restriction of $x$ to $J_{k}=(t_{k},t_{k+1}]$,
$k=0,1,...,m$. Set $\|\cdot\|_{b}$ be a seminorm in
$\mathrm{BM}_{h}'$ defined by
$$
\|x\|_{b}=\|x_{0}\|_{\mathrm{BM}_{h}}+\sup\left\{\|x(s)\|_{X},\,\
s\in[0,b]\right\}, \,\ x\in\mathrm{BM}'_{h}.
$$
\begin{lemma} (\cite{chang})\label{encadre}
Assume $x\in\mathrm{BM}'_{h}$, then for $t\in J$,
$x_{t}\in\mathrm{BM}_{h}$. Moreover,
$$
l\|x(t)\|_{X}\leq \|x_{t}\|_{\mathrm{BM}_{h}}\leq
\|x_{0}\|_{\mathrm{BM}_{h}}+l\sup\limits_{s\in[0,t]}\|x(s)\|,\,\
\text{where}\,\ l=\int_{-\infty}^{0}h(t)dt<+\infty.
$$
\end{lemma}

Before proceeding to the our main result, we shall set forth some
definitions and hypotheses from (\cite{pazy}, \cite{kellerman},
\cite{yosida}).
\subsection{Integrated Semigroups}

\begin{definition} (semigroup \cite{engel})
A semigroup is a set $S$ coupled with a binary operation $\ast$
($\ast: S\times S\longrightarrow S$) which is associative.
Associativity can also be realized as $F(F(x,y),z)=F(x,F(y,z))$
where $F(x,y)$ serves as the mapping from $S\times S$ to $S$.\\
A semigroup, unlike a group, need not have an identity element $e$
such that $x\ast e=x$, $\forall x\in S$.
\end{definition}
\begin{definition} ($C_{0}$ semigroup \cite{abada1}) A $C_{0}$ semigroup (or
strongly continuous semigroups) is a family,
$\{S(t)\}_{t\geq0}=\{S(t) \,\ | \,\ t\in \mathbb{R}_{+}\}$, of
bounded linear
operators from $X$ to $X$ satisfying:\\

(i) $S(0)=0$;\\

(ii) $t\longrightarrow S(t)$ is strongly continuous
($\lim\limits_{t\longrightarrow 0^{+}}S(t)f=f$ for each $f\in X$
with respect to the norm on $X$);

(iii) $S(s+t)=S(s)S(t) = \int_{0}^{s}(S(t+\sigma)-S(\sigma))d\sigma$
, for all $t,s \geq0.$
\end{definition}
\begin{definition} (Generator \cite{belleni}) Let $S$ be a semigroup. The (infinitesimal) generator of $S$, denoted by $A$, is given by the
equation:
$$
Af=\lim\limits_{t\longrightarrow 0^{+}} A_{t}f=
\lim\limits_{t\longrightarrow 0^{+}}\frac{S(t)f-f}{t},
$$
where the limit is evaluated in terms of the norm on $X$ and $f$ is
in the domain of $A$ if this limit exists.
\end{definition}

\begin{definition}
Let $\rho(A)$ is the resolvent set of an operator $A$ and $I$ is the
identity operator in $X$. Then $A$ is called a generator of an
integrated semigroup, if there exists $\delta\in\mathbb{R}$ such
that $(\delta,+\infty)\subset\rho(A)$, and there exists a strongly
continuous exponentially bounded family $\{S(t)\}_{t\geq0}$ of
linear bounded operators such that $S(0)=0$ and
$R(\lambda,A):=(\lambda
I-A)^{-1}=\lambda\int_{0}^{+\infty}\exp(-\lambda t)S(t)dt$ exists
for all $\lambda>\delta$.
\end{definition}
\begin{definition}
An integrated semigroup $\{S(t)\}_{t\geq0}$ is called locally
Lipschitz continuous if, for all $\delta>0$ there exists a constant
$\gamma$ such that
$$
\|S(t)-S(s)\|_{X}\leq\gamma|t-s|, \,\ t,s\in[0,\delta].
$$
Moreover it is called non-degenerate if $S(t)x=0$, for all $t\geq0$,
implies $x=0$.
\end{definition}
\begin{definition}\label{hy}
We say that a linear operator $A$ satisfies the Hille-Yosida
condition if there exists $\overline{M}\geq1$ and
$\delta\in\mathbb{R}$ such that $(\delta,+\infty)\subset\rho(A)$ and
$$
\sup\left\{(\lambda-\delta)^{n}\|R(\lambda,A)^{n}\|_{X}:
n\in\mathbb{N}, \lambda>\delta\right\}\leq \overline{M}.
$$
\end{definition}
\begin{theorem}(\cite{kellerman}, p.166) The following assertions are equivalent:\\

(i) $A$ is the generator of a non-degenerate, locally Lipschitz
continuous integrated semi-group;\\

(ii) $A$ satisfies the Hille-Yosida condition.
\end{theorem}

Throughout this paper we assume that:

The operator $A$ satisfies the Hille-Yosida condition: ($H_{0}$).

\begin{lemma}\label{satisf}
 The Operator $A$ satisfies the Hille-Yosida condition on $X$
(with $M=1$ and $\delta=0$). Therefore it follows that the operator
$A$ generates an integrated semigroup $S(t)$, $t\leq b$. Such that
$\|S'(t)\|\leq e^{-\kappa t}$, for $t\geq0$ and some constant
$\kappa>0$.
\end{lemma}

From \cite{kellerman}, since the operator $A$ hold the condition in
Definition \ref{hy}, this operator is the generator of a locally
Lipshitz continuous integrated semigroup $\{S(t)\}_{t\geq0}$
generates $C_{0}$ semigroup (see \cite{pazy}) on $\overline{D(A)}$
such that $$\|S(t)x\|_{X}\leq \overline{M}e^{\delta t}\|x\|_{X},$$
for all $t\geq0$ and $x\in \overline{D(A)}$. Let $A_{0}$ be the
generator of $C_{0}$ semigroup $\{S'(t)\}_{t\geq0}$ then $A_{0}$ is
the part of $A$ in $\overline{D(A)}$ defined by
$$
\left\{
\begin{array}{ll}
D(A_{0})=\{x\in D(A): Ax\in\overline{D(A)}\},\\
A_{0}x=Ax, \,\ \text{for}\,\ x\in \overline{D(A_{0})},
\end{array}
\right.
$$
we also have $\|S'(t)\|_{X}\leq \overline{M} e^{\delta t}$,
$t\geq0$, where $\overline{M}$ and $\delta$ are the constants
considered in Definition \ref{hy}.

\begin{lemma} (\cite{pazy})
Let $S(t)$ be a uniformly continuous semigroup of bounded linear
operators. Then $t\longrightarrow S(t)$ is differentiable in norm
and
$$
S'(t)=AS(t)=S(t)A.
$$
\end{lemma}

\section{Pseudo almost periodicity, almost periodicity, preliminaries and hypotheses}
It follows that the solution $x(\cdot)$ of
(\ref{elsevier1})-(\ref{elsevier111}) is from the space
$\mathrm{BM_{h}}$ and thus one may adopt and introduce the concept
of almost periodicity for sequences and continuous functions.
\begin{definition} (\cite{amerio}, \cite{diagana})
Let $f\in BC(\mathbb{R},X).$ We say that $f$ is almost periodic or
uniformly almost periodic (u.a.p), when the following property is
satisfied:
$$
\forall \varepsilon >0,\exists l_{\varepsilon }>0,\forall \alpha \in
\mathbb{ R},\exists \delta \in \left[ \alpha ,\alpha +l_{\varepsilon
}\right) ,\left\Vert f(\cdot +\delta )-f(\cdot )\right\Vert _{X
}<\varepsilon.
$$
\end{definition}
\begin{definition} Let $W$ a Banach space, and $\Omega$ a non empty subset
of the Banach space $Y\times Z$. We denote by $P_{c}(X)$ the set of
the compact subset of $Y\times Z$. Let a mapping $F\in C^{0}(J\times
Y\times Z,W)$; $F:(t,x)\longrightarrow F(t,x)$. We say that $F$ is
uniformly almost periodic in $t$ for $(y,z)\in \Omega$ when
$$
\forall \epsilon>0, \forall K\in P_{c}(Y\times Z), \exists
l=l(\epsilon,K),\,\ \forall \alpha\in\mathbb{R}, \,\ \exists
\tau\in[\alpha,\alpha+l),
$$
such that
$$
\sup\limits_{t\in\mathbb{R}}\sup\limits_{(y,z)\in
K}\|F(t+\tau,y,z)-F(t,y,z)\|\leq\epsilon.
$$
\end{definition}

\begin{definition} (\cite{zabut})
The set of sequences $\{\tau_{k+j}-\tau_{k}\}$, $k,j\in\mathbb{Z}$
is said to be uniformly almost periodic if for arbitrary
$\epsilon>0$ there exists a relatively dense set of
$\epsilon$-almost periods common for any sequences.
\end{definition}

\begin{definition} (\cite{farouk})
A sequence $I:\mathbb{Z}\longrightarrow X$ is called almost periodic
sequence, when for each $\epsilon>0$, $\exists
N_{\epsilon}\in\mathbb{N}$ such that among any $N_{\epsilon}$
consecutive integers there exists an integer $\tau$ with the
property
$$
\forall n\in\mathbb{Z}, \,\ \|I_{n+\tau}-I_{n}\|_{X}\leq\epsilon.
$$
In other words, $\widehat{T}(I,\epsilon)=\{\tau\in\mathbb{Z}:\,\
\forall n\in\mathbb{Z}\,\ \|I_{n+\tau}-I_{n}\|_{X}\leq\epsilon\}$ is
relatively dense in $\mathbb{Z}$. The class of all almost periodic
sequences is denoted by $APS$.
\end{definition}

\begin{definition} (\cite{zabut})
The function $x\in \mathrm{BM}'_{h}$ is
said to be piecewise almost periodic if \\

$(a_{1})$ the set of sequences $\{\tau_{k+j}-\tau_{k}\}$,
$k,j\in\mathbb{Z}$, $\left\{\tau_{k}\right\}\in \mathbf{T}$ is
uniformly almost periodic.\\

$(a_{2})$ for any $\epsilon>0$ there exists a real number $\delta>0$
such that if the points $t'$ and  $t''$ belong to one and the same
interval of continuity of $x(\cdot)$ and satisfy the inequality
$\left|t'-t''\right|<\delta$, then
$\left|x(t')-x(t'')\right|<\epsilon$.\\

$(a_{3})$ for any $\epsilon>0$ there exists a relatively dense set
$\mathrm{K}$ such that if $\tau\in\mathrm{K}$, then
$\left|x(t+\tau)-x(t)\right|<\epsilon$ for all $t\in\mathbb{R}$
satisfying $\left|t-t_{k}\right|>\epsilon$, $k\in\mathbb{Z}$.
\end{definition}

In other words, a Bohr almost periodic function is a continuous
function which possesses very much almost periods i.e. for all
$\epsilon>0$ the set
$$
T(f,\epsilon)=\{\tau\in \mathbb{Z}:
\|f(\cdot+\tau)-f(\cdot)\|_{X}\leq\epsilon\}
$$
is relatively dense in $\mathbb{Z}$. The elements of $T(f,\epsilon)$
are called $\epsilon$-periods. We denote by $AP(\mathbb{R},X)$ the
set of the Bohr a.p. functions from $\mathbb{R}$ to $X$. It is
well-known that the set $AP(\mathbb{R},X)$ is a Banach space with
the supremum norm. We refer the reader to (\cite{amerio},
\cite{cherif}, \cite{fink}) for the basic theory of almost periodic
functions and their applications. Besides, the concept of pseudo
almost periodicity (pap) was introduced by Zhang (see for example
\cite{zhang}) in the early nineties. It is a natural generalization
of the classical almost periodicity. Define the class of functions
$PAP_{0}(\mathbb{R},X)$ as follows.
$$
PAP_{0}(\mathbb{R},X)=\left\{f\in BC(\mathbb{R},X)\,\ \mid \,\
\lim\limits_{T\longrightarrow+\infty}
\frac{1}{2T}\int_{-T}^{T}\|f(t)\|_{X}dt=0 \right\}.
$$
\begin{definition}
A function $f\in BC(\mathbb{R},X)$ is called pseudo almost periodic
if it can be expressed as
$$
f=f_{1}+f_{2},
$$
where $f_{1}\in AP(\mathbb{R},X)$ ($AP(\mathbb{R}\times \Omega,X)$)
and $f_{2}\in PAP_{0}(\mathbb{R},X)$ ($PAP_{0}(\mathbb{R}\times
\Omega,X)$).
\end{definition}
The collection of such functions will be denoted by
$PAP(\mathbb{R},X)$ (resp $PAP(\mathbb{R}\times X,X)$. Besides, we
have a similar definition for pseudo almost periodic sequence.\\
\begin{definition}
A sequence $I:\mathbb{Z}\longrightarrow X$ is said to be in
$AP_{0}S$, when it satisfies the ergodicity condition
$$
\lim\limits_{p\longrightarrow+\infty}
\frac{1}{2p}\sum\limits_{i=-p}^{p}\|I(i)\|_{X}=0.
$$
\end{definition}
\begin{remark}
\end{remark} Notice that
\begin{enumerate}
  \item A sequence vanishing at infinity is a $AP_{0}S$ sequence.
  \item The sequence $(I(n))_{n\in\mathbb{Z}}$ defined by
\begin{equation*}I(n)= \left\{
\begin{array}{ll}
1, \,\ \text{if}\,\ n=2^{k},\\
0, \,\ \text{if}\,\ n\neq2^{k},
\end{array}\right.
\end{equation*} is an example of $AP_{0}S$ sequence which not
vanishing at infinity.
  \item For $k\in \mathbb{N}$, the sequence
  $(I(n))_{n\in\mathbb{Z}}$ defined by
  \begin{equation*}I(n)= \left\{
\begin{array}{ll}
1, \,\ \text{if}\,\ n=2^{k^{2}},\\
0, \,\ \text{if}\,\ n\neq2^{k^{2}},
\end{array}\right.
\end{equation*} is an example of an unbounded $AP_{0}S$ sequence.
\end{enumerate}
\begin{definition} We define the space of pseudo almost periodic sequences by
$$PAPS=APS\oplus AP_{0}S.$$
\end{definition}
\begin{remark}
If $f,g\in PAP(\mathbb{R}\times\Omega,\mathbb{R}^{n})$, then $f\pm
g\in PAP(\mathbb{R}\times\Omega,\mathbb{R}^{n})$ and $f\times g\in
PAP(\mathbb{R}\times\Omega,\mathbb{R}^{n})$.
\end{remark}
For more details about the properties of pseudo almost periodic
sequences we refer the reader to (\cite{farouk}, \cite{tran}).\\

Now, we introduce necessary fundamental properties of the almost
periodic concept, which will be used later.

\begin{lemma}\label{samoi} (\cite{samoil})
Let the set of sequences
$\left\{\tau_{j+k}-\tau_{k}\right\}_{k\in\mathbb{Z}}$,
$j\in\mathbb{Z}$ be uniformly almost periodic. Then for each $p>0$
there exists a positive integer $N$ such that on each interval of
length $p$ no more than $N$ elements of the sequence
$\left\{\tau_{k}\right\}$, i.e.,
$$
i(s,t)\leq N(t-s) +N,
$$
where $i(s,t)$ is the number of points $\tau_{k}$ in the interval
$(s,t)$
\end{lemma}
\begin{lemma}\label{para}\emph{(\cite{cherif})}\label{parametric}
Let the parametric function $t\longrightarrow\varphi(t,x(\cdot))$ is
almost periodic function, for all $x\in X\subseteq \mathbb{R}^{n}$.
Then the function $t\longrightarrow\int_{0}^{t}\varphi(t,x(s))ds$ is
almost periodic.
\end{lemma}

\begin{lemma} (\cite{alonso}, \cite{alonso1})\label{alalonso}
If $\{x(n)\}_{n\in\mathbb{Z}}$ is a $AP_{0}S$ sequence, then there
exists a function $w\in PAP_{0}(\mathbb{R},X)$ such that
$w(n)=x(n)$, $n\in\mathbb{Z}$.
\end{lemma}

For some preliminary results on almost periodic functions, we refer
the reader to (\cite{besicovitch}, \cite{fink} and \cite
{yoshizawa}).

\begin{definition} (See \cite{sadov})
An operator $F$ form Banach space $X$ to a Banach space $Y$ is
called a condensing operator if it is continuous and for every
bounded subset $B$ of $X$ the inequality $\chi[f(B)]<\chi(B)$ holds,
where $\chi(\cdot)$ denotes the measure of noncompactness.
\end{definition}
\begin{theorem} (Sadovskii \cite{sadov})
Let $F$ is a condensing operator on a Banach space $X$. If
$FD\subseteq D$ for a convex, bounded and closed subset $D$ of $X$,
then $F$ has a fixed point in $D$.
\end{theorem}
\begin{definition}
The function $x\in\mathrm{BM}'_{h}$ is said to be piecewise pseudo
almost periodic
solution of (\ref{elsevier1})-(\ref{elsevier111}) if $x$ is pseudo almost periodic function and satisfy;\\

(i) $x(t)=\phi(t), \,\ t\in(-\infty,0]$;\\

(ii)
$\left[x(s)-g\left(s,x_{s},\int_{0}^{s}h(s,\tau,x_{\tau})d\tau\right)\right]ds\in
D(A)$ on $J$;\\

(iii)
\begin{equation}\left\{
\begin{array}{ll}
x(t)=S'(t)\left[\phi(0)-g(0,\phi,0)\right]+g\left(t,x_{t},\int_{0}^{t}h(t,s,x_{s})ds\right)\\
+\lim\limits_{\lambda\longrightarrow+\infty}\int_{0}^{t}S'\left(t-s\right)\mathrm{B}_{\lambda}f\left(s,x_{s},\int_{0}^{s}k(t,\tau,x_{\tau})d\tau\right)ds\\
+\sum\limits_{0<t_{k}<t}S'\left(t-t_{k}\right)\delta(t-t_{k})I_{k}\left(x(t_{k}^{-})\right),
\,\ t\in J, \,\ \text{where}\,\ \mathrm{B}_{\lambda}=\lambda
R(\lambda,A).
\end{array}\right.
\end{equation}
\end{definition}

\begin{lemma} (\cite{samoil},
\cite{stamovalzabout})\label{dense} Suppose the following
conditions hold.\\

\emph{($A_{1}$)} The set of sequences $\{\tau_{k}^{j}\}$,
$\tau_{k}^{j}=\tau_{k+j}-\tau_{k}$, $k,j\in\mathbb{Z}$,
$\{\tau_{k}\}\in \mathbf{T}$, is almost periodic and there exists
$\theta>0$ such that
$$
\inf\limits_{k}\{\tau_{k}^{1}\}=\theta>0.
$$

\emph{($A_{2}$)} The functions $f$ and $g$ are almost periodic and
locally H\"older continuous with points of discontinuous at the
moments $t=\tau_{k}$, $k\in\mathbb{Z}$, at which it is continuous
from the
left.\\

\emph{($A_{3}$)} The sequence $b_{k}$, $k\in\mathbb{Z}$, is almost
periodic.\\
Then for each $\epsilon>0$ there exist $\epsilon_{1}$,
$0<\epsilon_{1}<\epsilon$ and relatively dense sets $T$ of real
numbers and $Q$ of whole numbers such that the following relations
fulfilled:\\

(i)$\|f(t+\tau,\cdot,\cdot)-f(t,\cdot,\cdot)\|_{X}< \epsilon, \,\ t\in\mathbb{R}, \,\ \tau\in T, \,\ \|t-\tau_{k}\|>\epsilon, \,\ k\in\mathbb{Z}$;\\

(ii)$\|g(t+\tau,\cdot,\cdot)-g(t,\cdot,\cdot)\|_{X}< \epsilon, \,\ t\in\mathbb{R}, \,\ \tau\in T, \,\ \|t-\tau_{k}\|>\epsilon, \,\ k\in\mathbb{Z}$;\\

(iii)$\|I_{k+q}-I_{k}\|_{X}< \epsilon, \,\ q\in Q, \,\ k\in\mathbb{Z}$;\\

(iv) $\|\tau_{k+q}-\tau_{k}\|< \epsilon_{1}, \,\ q\in Q, \,\ \tau\in
T, \,\ k\in\mathbb{Z}$.
\end{lemma}

The solutions of the functional impulsive system of
integro-differential equations of the form
(\ref{elsevier1})-(\ref{elsevier111}) are characterized in the
following way:

\begin{equation}\label{elsevier2}x(t)=\left\{
\begin{array}{ll}
\phi(t), \,\ \text{if}\,\ t\in(-\infty,0],\\
S'(t)\left[\phi(0)-g(0,\phi,0)\right]+g\left(t,x_{t},\int_{0}^{t}h(t,s,x_{s})ds\right)\\
+\lim\limits_{\lambda\longrightarrow+\infty}\int_{0}^{t}S'\left(t-s\right)\mathrm{B}_{\lambda}f\left(s,x_{s},\int_{0}^{s}k(t,\tau,x_{\tau})d\tau\right)ds\\
+\sum\limits_{0<t_{k}<t}S'\left(t-t_{k}\right)\delta(t-t_{k})I_{k}\left(x(t_{k}^{-})\right),
\,\ t\in J, \,\ \text{where}\,\ \mathrm{B}_{\lambda}=\lambda
R(\lambda,A),
\end{array}\right.
\end{equation}
where
\begin{equation*}\delta(t-t_{k})=\left\{
\begin{array}{ll}
1, \,\ \text{if}\,\ t=t_{k},\\
0, \,\ \text{if}\,\ t\neq t_{k}.\\
\end{array}\right.
\end{equation*}



Now, we are in a position to state the existence theorem. Our first
theorem is based on the Sadovskii's theorem.

\section{Existence of pseudo almost periodic solution of impulsive integro-differential systems with infinite
Delay}
\begin{lemma}\label{boch}
The functional $\xi: s\longrightarrow
S'\left(t-s\right)\mathrm{B}_{\lambda}f\left(s,x_{s},\int_{0}^{s}k(t,\tau,x_{\tau})d\tau\right)$
is integrable on $[0,t)$.
\end{lemma}
Proof.\\
We have, from the hypothesis ($H_{2}$):
\begin{eqnarray*}
\|\xi(s)\|&=&\left\|S'\left(t-s\right)\mathrm{B}_{\lambda}f(s,x_{s},\int_{0}^{s}k(t,\tau,x_{\tau})d\tau)\right\|\\
&=&\left\|AS\left(t-s\right)\mathrm{B}_{\lambda}f(s,x_{s},\int_{0}^{s}k(t,\tau,x_{\tau})d\tau)\right\|\\
&\leq&\frac{\overline{M}\lambda}{\lambda-\delta}
\left\|(-A)^{1-\beta}S(t-s)\right\|\left\{\left\|(-A)^{\beta}f\left(s,x_{s},\int_{0}^{s}k(t,\tau,x_{\tau})d\tau\right)
-(-A)^{\beta}f(s,0,0)\right\|\right.\\
&+&\left.c_{2}\right\}\\
&\leq&\frac{\overline{M}\lambda}{\lambda-\delta}
\frac{C_{1-\beta}}{(t-s)^{1-\beta}}\left\{F_{1}\|x_{s}\|_{\mathrm{BM}_{h}}\right.\\
&+&\left.c_{2}\right\}\\
&\leq&\frac{\overline{M}\lambda}{\lambda-\delta}
\frac{C_{1-\beta}}{(t-s)^{1-\beta}}\left\{F_{1}(\|x_{0}\|_{\mathrm{BM}_{h}}+l\sup\limits_{r\in[0,s]}\|x(r)\|_{X})
+c_{2} \right\}.
\end{eqnarray*}

Thus from Bochner's theorem, it follows that $s\longrightarrow
S'\left(t-s\right)\mathrm{B}_{\lambda}f\left(s,x_{s},\int_{0}^{s}k(t,\tau,x_{\tau})d\tau\right)$
is integrable on $[0,t)$. This achieve the proof of Lemma
\ref{boch}.

 Define the operator
$\Gamma:\mathrm{BM}''_{h}\longrightarrow\mathrm{BM}''_{h}$ by
\begin{equation*}\Gamma(x)(t)=\left\{
\begin{array}{ll}
\phi(t), \,\ \text{if}\,\ t\in(-\infty,0],\\
-S'(t)[\phi(0)-g(0,\phi,0)]+g\left(t,x_{t},\int_{0}^{t}h(t,s,x_{s})ds\right)\\
+\lim\limits_{\lambda\longrightarrow+\infty}\int_{0}^{t}S'\left(t-s\right)\mathrm{B}_{\lambda}f(s,x_{s},\int_{0}^{s}k(t,\tau,x_{\tau})d\tau)ds\\
+\sum\limits_{0<t_{k}<t}S'\left(t-t_{k}\right)\delta(t-t_{k})I_{k}\left(x(t_{k}^{-})\right),
\,\ t\in J,
\end{array}\right.
\end{equation*}
The operator $\Gamma$ has a fixed point $x(\cdot)$. This fixed point
is then PPAP solution of the system
(\ref{elsevier1})-(\ref{elsevier111}).\\

We show the following

\begin{lemma}\label{1}
For $u$ pseudo almost periodic the operator $\Gamma(u)$ is pseudo
almost periodic.
\end{lemma}
Proof.\\
First, since $\alpha_{k}\in L^{1}(J,[0,+\infty))$, by hypothesis
($H_{5}$) using the Hille-Yosida condition with $n=1$, (we have, $
\left\|B_{\lambda}\right\|_{X}=\left\|\lambda(\lambda
I-A)^{-1}\right\|_{X}\leq
\frac{\overline{M}\lambda}{\lambda-\delta}$), and by hypothesis
($H_{2}$)  we observe that $\Gamma$ is bounded.
\begin{eqnarray*}
\|\Gamma (x)(t)\|_{X}&\leq& \left\|S'(t)\right\|_{X}\|\phi(0)-g(0,\phi,0)\|_{X}+\left\|g\left(t,x_{t},\int_{0}^{t}h(t,s,x_{s})ds\right)\right\|_{X}\\
&+&
\lim\limits_{\lambda\longrightarrow+\infty}\frac{\overline{M}\lambda}{\lambda-\delta}\left\|\int_{0}^{t}S'\left(t-s\right)f\left(s,x_{s},\int_{0}^{s}k(t,\tau,x_{\tau})d\tau\right)ds\right\|_{X}\\
&+&
\left\|\sum\limits_{0<t_{k}<t}S'\left(t-t_{k}\right)I_{k}\left(x(t_{k}^{-})\right)\right\|_{X}
\\
&\leq& M
\|\phi(0)-g(0,\phi,0)\|_{X}+\left\|g\left(t,x_{t},\int_{0}^{t}h(t,s,x_{s})ds\right)-g(0,0,0)\right\|_{X}+\|g(0,0,0)\|_{X}\\
&+& \lim\limits_{\lambda\longrightarrow+\infty}\frac{M\overline{M}\lambda}{\lambda-\delta}\left\|\int_{0}^{t}f\left(s,x_{s},\int_{0}^{s}k(t,\tau,x_{\tau})d\tau\right)ds\right\|_{X}\\
&+&
\left\|\sum\limits_{0<t_{k}<t}S'\left(t-t_{k}\right)I_{k}\left(x(t_{k}^{-})\right)\right\|_{X}\\
&\leq& M
\|\phi(0)-g(0,\phi,0)\|_{X}+G_{1}\|x_{t}\|_{\mathrm{BM}_{h}}+G_{2}p(t)\|x_{t}\|_{\mathrm{BM}_{h}}+\|g(0,0,0)\|_{X}\\
&+& M\overline{M}\int_{0}^{t}\sup\limits_{\|x_{s}\|_{C},\left\|\int_{0}^{s}k(t,\tau,x_{\tau})d\tau\right\|<k}\left\|f\left(s,x_{s},\int_{0}^{s}k(t,\tau,x_{\tau})d\tau\right)\right\|_{X} ds\\
&+&
\left\|\sum\limits_{0<t_{k}<t}S'\left(t-t_{k}\right)I_{k}\left(x(t_{k}^{-})\right)\right\|_{X}
\\ &\leq&
\|\phi(0)-g(0,\phi,0)\|+(G_{2}\bar{p}+G_{1})\left(\|x_{0}\|_{\mathrm{BM}_{h}}+\sup\limits_{s\in[0,b]}\|x(s)\|\right)+\|g(0,0,0)\|_{X}\\
&+& M\overline{M}\int_{0}^{b}\alpha_{k}(s)ds+
M\left\|\sum\limits_{0<t_{k}<b}L_{k}\left(x(t_{k})\right)\right\|_{X}\\
&<& \infty
\end{eqnarray*}
Hence $\Gamma(x)$ is bounded.\\

Now we show that $\Gamma(x)(t)$ is pseudo almost periodic with
respect to $t\in\mathbb{R}$. Note that the functions $f,g,h$ and $k$
are pseudo almost periodic. By the fact that the sum and product of
two functions pseudo almost periodic are pseudo almost periodic and
by Lemma \ref{para}, we have
\begin{eqnarray*}
t&\longrightarrow& -S'(t)[\phi(0)-g(0,\phi,0)]+g\left(t,x_{t},\int_{0}^{t}h(t,s,x_{s})ds\right)\\
&+&\lim\limits_{\lambda\longrightarrow+\infty}\int_{0}^{t}S'\left(t-s\right)\mathrm{B}_{\lambda}f\left(s,x_{s},\int_{0}^{s}k(t,\tau,x_{\tau})d\tau\right)ds,
\end{eqnarray*}
is pseudo almost periodic.\\

Let us prove that the function
$$
z:t\longrightarrow
\sum\limits_{0<t_{k}<t}S'\left(t-t_{k}\right)\delta(t-t_{k})I_{k}\left(x(t_{k}^{-})\right),
$$
belongs to $PAP(\mathbb{R},X)$. Because $I_{k}$ is pseudo almost
periodic sequence, one have
$I_{k}(x(t_{k}^{-}))=I_{k}=I_{k}^{1}+I_{k}^{2}$ where $I_{k}^{1}\in
AP(\mathbb{R},X)$ and $I_{k}^{2}\in AP_{0}S(\mathbb{R},X)$.
Consequently,

\begin{eqnarray*}
z(t)&=&
\sum\limits_{0<t_{k}<t}S'\left(t-t_{k}\right)\delta(t-t_{k})I_{k}^{1}+\sum\limits_{0<t_{k}<t}S'\left(t-t_{k}\right)\delta(t-t_{k})I_{k}^{2}\\
&=& z_{1}(t)+z_{2}(t).
\end{eqnarray*}
Since $I_{k}^{1}$ is almost periodic, then for a fixed $\epsilon>0$
the set
$$
\widehat{T}(I_{k}^{1},\epsilon)=\left\{\tau\in\mathbb{R}:\,\
\|I_{k}^{1}(t+\tau)-I_{k}^{1}(t)\|_{X}\leq\epsilon \,\ \forall
n\in\mathbb{Z} \right\}
$$
is relatively dense in $\mathbb{Z}$. That is there exists a positive
integer $l_{\epsilon}$ such that any interval with the length
$l_{\epsilon}$ contains at least one point of
$\widehat{T}(I_{k}^{1},\epsilon)$. Now, let $\epsilon>0$, $\tau\in
T$, $q\in Q$ where the $T$ and $Q$ are defined as in Lemma
\ref{dense}. Then
\begin{eqnarray*}
\|z_{1}(t+\tau)-z_{1}(t)\|_{\alpha}&=&\left\|(-A)^{-\alpha}(z_{1}(t+\tau)-z_{1}(t))\right\|_{X}
\\
&=&\left\|(-A)^{-\alpha}\left(\sum\limits_{0<t_{k}<t+\tau}S'(t+\tau-t_{k})\delta(t-t_{k})I_{k}^{1}\right.\right.\\
&-&\left.\left.\sum\limits_{0<t_{k}<t}S'(t-t_{k})\delta(t-t_{k})I_{k}^{1}\right)\right\|_{X}\\
&\leq&
\sum\limits_{0<t_{k}<t}\left\|(-A)^{-\alpha}S'(t-t_{k})\right\|_{X}\left\|I_{k+q}^{1}-I_{k}^{1}\right\|_{X}\\
&\leq&
\sum\limits_{0<t_{k}<t}\left\|(-A)^{1-\alpha}S(t-t_{k})\right\|_{X}\left\|I_{k+q}^{1}-I_{k}^{1}\right\|_{X}
\\
&\leq&
\sum\limits_{0<t_{k}<t}C_{1-\alpha}(t-t_{k})^{1-\alpha}\exp(-\lambda(t-t_{k}))\left\|I_{k+q}^{1}-I_{k}^{1}\right\|_{X}\\
&\leq& \epsilon
\sum\limits_{0<t_{k}<t}C_{1-\alpha}(t-t_{k})^{1-\alpha}\exp(-\lambda(t-t_{k}))
\\
&\leq&
\epsilon\left[\sum\limits_{0<t-t_{k}<1}C_{1-\alpha}(t-t_{k})^{1-\alpha}\exp(-\lambda(t-t_{k}))\right.
\\
&+&\left.\sum\limits_{j=1}^{\infty}\sum\limits_{j<t-t_{k}<j+1}C_{1-\alpha}(t-t_{k})^{1-\alpha}\exp(-\lambda(t-t_{k}))\right]
\\
&\leq& \epsilon C_{1-\alpha} N
\left(\frac{m^{1-\alpha}}{\exp(-\lambda)}+\frac{1}{\exp(\lambda)-1}\right),
\end{eqnarray*}
 where $\|t-t_{k}\|>\epsilon$ and $m=\min\{t-t_{k},\,\
0<t-t_{k}\leq1\}$.\\

Now, we shall prove that
$$
\lim\limits_{r\longrightarrow+\infty}\frac{1}{2r}\int_{-r}^{r}\left\|z_{2}(s)\right\|_{\alpha}ds=0.
$$
By Lemma \ref{alalonso} there exists $w\in PAP_{0}(\mathbb{R},X)$
such that $w(k)=I^{2}(k)=I_{k}^{2}$, $k\in\mathbb{Z}$.
\begin{eqnarray*}
\lim\limits_{r\longrightarrow+\infty}\frac{1}{2r}\int_{-r}^{r}\left\|z_{2}(s)\right\|_{\alpha}ds
&=&\lim\limits_{r\longrightarrow+\infty}\frac{1}{2r}\int_{-r}^{r}\left\|\sum\limits_{0<t_{k}<t}S'(t-t_{k})I_{k}^{2}\right\|_{\alpha}ds\\
&=&\lim\limits_{r\longrightarrow+\infty}\frac{1}{2r}\int_{-r}^{r}\left\|\sum\limits_{0<t_{k}<t}(-A)^{-\alpha}S'(t-t_{k})I_{k}^{2}\right\|_{\alpha}ds
\\
&\leq&\lim\limits_{r\longrightarrow+\infty}\frac{1}{2r}\int_{-r}^{r}\left\|\sum\limits_{0<t_{k}<t}(-A)^{1-\alpha}S(t-t_{k})\right\|_{X}\left\|I_{k}^{2}\right\|_{X}ds\\
&\leq&\lim\limits_{r\longrightarrow+\infty}\frac{1}{2r}\int_{-r}^{r}\left\|\sum\limits_{0<t_{k}<t}C_{1-\alpha}(t-t_{k})^{1-\alpha}\exp(-\lambda(t-t_{k}))\right\|_{X}\left\|I_{k}^{2}\right\|_{X}ds
\\
&\leq&\lim\limits_{r\longrightarrow+\infty}\frac{C_{1-\alpha}}{2r}\int_{-r}^{r}\sum\limits_{0<t-t_{k}\leq1}(t-t_{k})^{1-\alpha}\exp(-\lambda(t-t_{k}))
\\
&&\times\left\|w(k)\right\|_{X}ds
\\
&+&\lim\limits_{r\longrightarrow+\infty}\frac{C_{1-\alpha}}{2r)}\int_{-r}^{r}\sum\limits_{j=1}^{+\infty}\sum\limits_{j<t-t_{k}\leq
j+1}(t-t_{k})^{1-\alpha}\exp(-\lambda(t-t_{k}))\\
&&\times\left\|w(k)\right\|_{X}ds
\\
&\leq&\lim\limits_{r\longrightarrow+\infty}\frac{C_{1-\alpha}}{2r}\int_{-r}^{r}\sum\limits_{0<t-t_{k}\leq1}(t-t_{k})^{1-\alpha}\exp(-\lambda(t-t_{k}))\\
&&\times\|w(t)I(t-k)\|_{X}ds
\\
&+&\lim\limits_{r\longrightarrow+\infty}\frac{C_{1-\alpha}}{2r}\int_{-r}^{r}\sum\limits_{j=1}^{+\infty}\sum\limits_{j<t-t_{k}\leq
j+1}(t-t_{k})^{1-\alpha}\exp(-\lambda(t-t_{k}))\\
&&\times\|w(t)I(t-k)\|_{X}ds
\\
&\leq&\lim\limits_{r\longrightarrow+\infty}\frac{C_{1-\alpha}}{2r}\int_{-r}^{r}\sum\limits_{0<t-t_{k}\leq1}(t-t_{k})^{1-\alpha}\exp(-\lambda(t-t_{k}))\|w(t)\|_{X}ds
\\
&+&\lim\limits_{r\longrightarrow+\infty}\frac{C_{1-\alpha}}{2r}\int_{-r}^{r}\sum\limits_{j=1}^{+\infty}\sum\limits_{j<t-t_{k}\leq
j+1}(t-t_{k})^{1-\alpha}\exp(-\lambda(t-t_{k})) \|w(t)\|_{X}ds
\end{eqnarray*}
\begin{eqnarray*}
&\leq&\lim\limits_{r\longrightarrow+\infty}\frac{C_{1-\alpha}}{2r}\int_{-r}^{r}N\frac{m^{1-\alpha}}{\exp(-\lambda)}\|w(t)\|_{X}ds
\\
&+&
\lim\limits_{r\longrightarrow+\infty}\frac{C_{1-\alpha}}{2r}\int_{-r}^{r}\frac{N}{\exp(\lambda)-1}\|w(t)\|_{X}ds=0.
\end{eqnarray*}

This achieve the proof of Lemma \ref{1}.\\

Now, for $\phi\in \mathrm{BM}_{h}$ define $\tilde{\phi}$ by
\begin{equation*}\tilde{\phi}(t)=\left\{
\begin{array}{ll}
\phi(t), \,\ \text{if}\,\ t\in(-\infty,0],\\
S'(t)\phi(0), \,\ \text{if}\,\ t\in J=[0,b],
\end{array}\right.
\end{equation*}
then $\tilde{\phi}\in\mathrm{BM}'_{h}$. Let
$x(t)=y(t)+\tilde{\phi}(t)$, $-\infty<t\leq b$. Observe that
$x(\cdot)$ satisfies (\ref{elsevier2}) if and only if $y$ satisfies
$y_{0}=0$ and
\begin{eqnarray*}
y(t)&=&
-S'(t)g(0,\phi,0)+g\left(t,y_{t},\int_{0}^{t}h(t,s,y_{s})ds\right)\\
&+&\lim\limits_{\lambda\longrightarrow+\infty}\int_{0}^{t}S'\left(t-s\right)\mathrm{B}_{\lambda}f(s,y_{s},\int_{0}^{s}k(t,\tau,y_{\tau})d\tau)ds\\
&+&\sum\limits_{0<t_{k}<t}S'\left(t-t_{k}\right)\delta(t-t_{k})I_{k}\left(y(t_{k}^{-})\right),
\,\ t\in J,
\end{eqnarray*}

Let $\mathrm{BM}''_{h}=\{y\in \mathrm{BM}'_{h}: y_{0}=0\in
\mathrm{BM}_{h}\}$, then for any $y\in\mathrm{BM}''_{h}$, we have
$$
\|y\|_{b}=\|y_{0}\|_{\mathrm{BM}_{h}}+\sup\{\|y(s)\|_{X}: 0\leq
s\leq b\},
$$
thus $(\mathrm{BM}''_{h},\|\cdot\|_{b})$ is a Banach space. Let
$q>0$ and define
$$
B_{q}=\{y\in\mathrm{BM}''_{h}: \|y\|_{b} \leq q\}.
$$
 Then for each $q>0$ the set $B_{q}$ is clearly bounded,
closed convex set in $\mathrm{BM}''_{h}$. Further, using the Lemma
\ref{encadre}, for any $y\in B_{q}$ we have
\begin{eqnarray}\label{adoc}
\|y_{t}+\tilde{\phi}_{t}\|_{\mathrm{BM}_{h}}&\leq&
\|y_{t}\|_{\mathrm{BM}_{h}}+\|\tilde{\phi}_{t}\|_{\mathrm{BM}_{h}}\nonumber\\
&\leq& \|y_{0}\|_{\mathrm{BM}_{h}}+l
\sup_{s\in[0,t]}\|(s)\|+\|\tilde{\phi}_{0}\|_{\mathrm{BM}_{h}}+l
\sup_{s\in[0,t]}\|\tilde{\phi}(s)\|\nonumber\\
&\leq& l
q+\|\tilde{\phi}_{0}\|_{\mathrm{BM}_{h}}+l\sup\limits_{s\in[0,b]}\|T(s)\phi(0)\|\nonumber\\
&\leq& l(q+M\|\phi(0)\|_{X})+\|\phi\|_{\mathrm{BM}_{h}}=q'.
\end{eqnarray}
Further by using Lemma \ref{encadre} and (\ref{adoc}), for each
$t\in J$ we have
$$
\|y(t)+\tilde{\phi}(t)\|_{X}\leq
l^{-1}\|y_{t}-\tilde{\phi}_{t}\|_{\mathrm{BM}_{h}}\leq l^{-1}q'.
$$
This give,
\begin{equation}\label{used}
\sup\limits_{t\in J}\|y(t)+\tilde{\phi}(t)\|\leq l^{-1}q'.
\end{equation}
By using hypothesis ($H_{5}$) and (\ref{used}), for each
($k=1,...,m$) we obtain
\begin{equation}\label{eqjump}
\|I_{k}(y(t_{k}^{-})+\tilde{\phi}(t_{k}^{-}))\|_{X}\leq
L_{k}(\|y(t_{k}^{-})+\tilde{\phi}(t_{k}^{-})\|_{X})\leq
L_{k}\left(\sup\limits_{t\in J}\|y(t)+\tilde{\phi}(t)\|\right)\leq
L_{k}(l^{-1}q').
\end{equation}
Define the operator $\psi:\mathrm{BM}''_{h}\longrightarrow
\mathrm{BM}''_{h}$ by
\begin{equation*}\psi(y)(t)=\left\{
\begin{array}{ll}
0, \,\ \text{if}\,\ t\in(-\infty,0],\\
-S'(t)g(0,\phi,0)+g\left(t,y_{t},\int_{0}^{t}h(t,s,y_{s})ds\right)\\
+\lim\limits_{\lambda\longrightarrow+\infty}\int_{0}^{t}S'\left(t-s\right)\mathrm{B}_{\lambda}f(s,y_{s},\int_{0}^{s}k(t,\tau,x_{\tau})d\tau)ds\\
+\sum\limits_{0<t_{k}<t}S'\left(t-t_{k}\right)\delta(t-t_{k})I_{k}\left(x(t_{k}^{-})\right),
\,\ t\in J,
\end{array}\right.
\end{equation*}
it's clear that the operator $\Gamma$ has a fixed point if and only
if $\psi$ has a fixed point. Thus, we prove that $\psi$ has a fixed
point.\\

Now, we prove some necessary lemmas:

\begin{lemma}
$\psi(B_{q})\subseteq B_{q}$ for some $q>0$.
\end{lemma}
Proof.\\
We claim that there exists a positive integer $q$ such that
$\psi(B_{q})\subseteq B_{q}$. If this is not true, then for each
positive integer $q$, there is a function $y^{q}\in B_{q}$, but
$\psi(y^{q})\notin B_{q}$, that is, $\|\psi(y^{q}(t))\|_{X}>q$ for
some $t(q)\in J$, where $t(q)$ denotes $t$ depending on $q$.
However, on the other hand, using ($H_{1}$)-($H_{4}$) and the
inequality (\ref{eqjump}) we have,
\begin{eqnarray}\label{qu}
q&<& \left\|\psi(y^{q})(t)\right\|_{X}\nonumber\\
&\leq&
\left\|S'(t)\right\|_{X}\left(\|g(0,\phi,0)-g(t,0,0)\|_{X}+\|g(t,0,0)\|_{X}\right)\nonumber\\
&+&
\left\|g\left(t,y_{t}^{q}+\widetilde{\phi}_{t},\int_{0}^{t}h(t,s,y_{s}^{q}+\widetilde{\phi}_{s})ds\right)-g(t,0,0)\right\|_{X}+\left\|g(t,0,0)\right\|_{X}\nonumber\\
&+&
\left\|\lim\limits_{\lambda\longrightarrow+\infty}\int_{0}^{t}S'(t-s)B_{\lambda}f\left(s,y_{s}^{q}+\widetilde{\phi}_{s},\int_{0}^{s}k(t,\tau,y_{\tau}^{q}+\widetilde{\phi}_{\tau})d\tau\right)ds\right\|_{X}\nonumber\\
&+& \sum\limits_{0<
t_{k}<t}\left\|S'(t-t_{k})\right\|_{X}\left\|I_{k}\left(y^{q}(t_{k}^{-})+\widetilde{\phi}(t_{k}^{-})\right)\right\|_{X}\nonumber\\
&\leq& M \left(G_{1}\|\phi\|_{\mathrm{BM}_{h}}+\sup\limits_{t\in
J}\|g(t,0,0)\|\right)+G_{1}\left\|y_{t}^{q}+\widetilde{\phi}_{t}\right\|_{\mathrm{BM}_{h}}\nonumber\\
&+&G_{2}\left\|\int_{0}^{t}h(t,s,y_{s}^{q}+\widetilde{\phi}_{s})ds\right\|_{X}+\sup\limits_{t\in
J} \|g(t,0,0)\|\nonumber\\
&+&
\left\|\lim\limits_{\lambda\longrightarrow+\infty}\int_{0}^{t}S'(t-s)B_{\lambda}f\left(s,y_{s}^{q}+\widetilde{\phi}_{s},\int_{0}^{s}k(t,\tau,y_{\tau}^{q}+\widetilde{\phi}_{\tau})d\tau\right)ds\right\|_{X}+M\sum\limits_{k=1}^{m}L_{k}(l^{-1}q').\nonumber
\end{eqnarray}
Using the Hille-Yosida condition with $n=1$, we have,
$$
\left\|B_{\lambda}\right\|_{X}=\left\|\lambda(\lambda
I-A)^{-1}\right\|_{X}\leq
\frac{\overline{M}\lambda}{\lambda-\delta}.
$$
Further $B_{\lambda}x\longrightarrow x$ as
$\lambda\longrightarrow+\infty$ for all $x\in \overline{D(A)}$,
\begin{equation}\label{mbar}
\lim\limits_{\lambda\longrightarrow+\infty}\|B_{\lambda}x\|_{X}\leq
\overline{M}\|x\|_{X} \,\ \text{and} \,\
\lim\limits_{\lambda\longrightarrow+\infty}\|B_{\lambda}\|_{X}\leq\overline{M}.
\end{equation}
From (\ref{mbar}) and hypothesis ($H_{3}$), for each $t\in J$, we
have
\begin{multline}\label{omega}
\sup\left\{\left\|y_{t}+\widetilde{\phi}_{t}\right\|_{\mathrm{BM}_{h}},\left\|\int_{0}^{t}h\left(t,s,y_{s}+\widetilde{\phi}_{s}\right)ds\right\|_{X},\left\|\int_{0}^{b}k\left(t,s,y_{s}+\widetilde{\phi}_{s}\right)ds\right\|_{X}\right\}\\
\leq q'\max\limits_{t\in J}\{1,p(t),q(t)\}=
q'\overline{\omega}.\phantom{+++++++++++}
\end{multline}
Taking $G_{3}=\sup\limits_{t\in J}\|g(t,0,0)\|$, using inequalities
(\ref{mbar}), (\ref{omega}) and the hypothesis ($H_{4}$)-($H_{5}$),
from (\ref{qu}), we obtain
$$
q<
M\left(G_{1}\|\phi\|_{\mathrm{BM}_{h}}+G_{3}\right)+(G_{1}+G_{2})q'\overline{\omega}+G_{3}+M\overline{M}\int_{0}^{b}\alpha_{q'\overline{\omega}}(s)ds+M\sum\limits_{k=1}^{m}L_{k}(l^{-1}q').
$$
Therefore
\begin{eqnarray}\label{qi2}
1&<&\frac{M(G_{1}\|\phi\|_{\mathrm{BM}_{h}}+G_{3})}{q}+(G_{1}+G_{2})\overline{\omega}\frac{q'}{q}+\frac{M\overline{M}q'\overline{\omega}}{q}\frac{1}{q'\overline{\omega}}\int_{0}^{b}\alpha_{q'\overline{\omega}}(s)ds\nonumber\\
&+&M\frac{l^{-1}q'}{q}\sum\limits_{k=1}^{m}\frac{L_{k}(l^{-1}q')}{l^{-1}q'}.
\end{eqnarray}
Noting that $\lim\limits_{q\longrightarrow+\infty}\frac{q'}{q}=l$,
therefore $q'\longrightarrow+\infty$ as $q\longrightarrow+\infty$.\\
Taking lower limit from (\ref{qi2}) and from hypothesis
($H_{4}$)-($H_{5}$) we get
$$
\left(G_{1}+G_{2}\right)\overline{\omega}l+M\overline{M}\overline{\omega}l\mu+M\sum\limits_{k=1}^{m}\lambda_{k}\geq1.
$$
This contradicts the hypothesis of Theorem \ref{th1}. Hence for some
positive integer $q$, $\psi(B_{q})\subseteq B_{q}$.\\

\begin{lemma} \emph{(Decomposition of $\psi$ \cite{sadov})}. The completely continuous operators,
contractions and sum of these two operators are condensing
operators. Thus to prove that the $\psi$ is condensing operator we
decompose $\psi$ as $\psi=\psi_{1}+\psi_{2}$, and prove that the
operator $\psi_{1}$ is a contraction while $\psi_{2}$ is a
completely continuous operator.\\
The operator $\psi_{1},\psi_{2}$ on $B_{q}$, respectively are
defined by
\begin{equation}
\psi_{1}y(t)=-S'(t)g(0,\phi,0)+g\left(t,y_{t}+\widetilde{\phi}_{t},\int_{0}^{t}h\left(t,s,y_{s}+\widetilde{\phi}_{s}\right)ds\right),
\,\ t\in J,
\end{equation}
\begin{eqnarray}\label{psi2}
\psi_{2}y(t)&=&
\lim\limits_{\lambda\longrightarrow+\infty}\int_{0}^{t}S'(t-s)B_{\lambda}f\left(s,y_{s}+\widetilde{\phi}_{s},\int_{0}^{s}k\left(t,\tau,y_{\tau}+\widetilde{\phi}_{\tau}\right)d\tau\right)ds\nonumber\\
&+&\sum\limits_{0<t_{k}<t}S'(t-t_{k})\delta(t-t_{k})I_{k}\left(y(t_{k}^{-})+\widetilde{\phi}(t_{k}^{-})\right),
\,\ t\in J.
\end{eqnarray}
\end{lemma}
\begin{lemma} $\psi_{1}$ is a contraction.
\end{lemma}
Proof.\\
Let any $u,v\in B_{q}$ and $t\in J$. Then using the hypotheses
($H_{1}$) and Lemma \ref{encadre}, we have
\begin{eqnarray*}
\left\|\psi_{1}u(t)-\psi_{1}v(t)\right\|_{X}&\leq&
\left\|g\left(t,u_{t}+\widetilde{\phi}_{t},\int_{0}^{t}h\left(t,s,u_{s}+\widetilde{\phi}_{s}\right)ds\right)\right.\\
&-&\left.g\left(t,v_{t}+\widetilde{\phi}_{t},\int_{0}^{t}h\left(t,s,v_{s}+\widetilde{\phi}_{s}\right)ds\right)\right\|_{X}\\
&\leq&
G_{1}\|u_{t}-v_{t}\|_{\mathrm{BM}_{h}}+G_{2}\int_{0}^{t}\left\|h\left(t,s,u_{s}+\widetilde{\phi}_{s}\right)-h\left(t,s,v_{s}+\widetilde{\phi}_{s}\right)\right\|_{X} ds\\
&\leq& \left(G_{1}+bG_{2}H\right)\|u_{t}-v_{t}\|_{\mathrm{BM}_{h}}\\
&\leq&
\left(G_{1}+bG_{2}H\right)\left(\|u_{0}\|_{\mathrm{BM}_{h}}+\|v_{0}\|_{\mathrm{BM}_{h}}+l\sup\limits_{s\in[0,t]}\|u(s)-v(s)\|\right).
\end{eqnarray*}
Since $\|u_{0}\|_{\mathrm{BM}_{h}}=\|v_{0}\|_{\mathrm{BM}_{h}}=0$,
we have
$$
\left\|\psi_{1}u(t)-\psi_{1}v(t)\right\|_{X}\leq\left(G_{1}+bG_{2}H\right)l\sup\limits_{s\in[0,t]}\|u(s)-v(s)\|\leq\left(G_{1}+bG_{2}H\right)l\|u-v\|_{b}.
$$
Therefore $\psi_{1}$ is a contraction on $B_{q}$.\\

\begin{lemma} $\psi_{2}$ maps $B_{q}$ in to an equicontinuous
family.
\end{lemma}
Proof.\\ Let any $y\in B_{q}$ and
$t_{1},t_{2}\in(-\infty,0]$. We prove the equicontinuity for these
case $0<t_{1}<t_{2}\leq b$, as the equicontinuity for these cases
$t_{2}<t_{1}\leq0$ and $t_{1}<0<t_{2}\leq b$ are obvious.

For $0<t_{1}<t_{2}\leq b$ and $\epsilon>0$ sufficiently small, using
the hypothesis ($H_{4}$) and the conditions (\ref{adoc}),
(\ref{eqjump}) and (\ref{omega}) we obtain
\begin{eqnarray*}
\left\|\psi_{2}(y)(t_{1})-\psi_{2}(y)(t_{2})\right\|_{X} &\leq&
\left\|\lim\limits_{\lambda\longrightarrow+\infty}\int_{0}^{t_{1}-\epsilon}\left(S'(t_{2}-s)-S'(t_{1}-s)\right)B_{\lambda}\right.\\
&\times&\left. f\left(s,y_{s}+\widetilde{\phi}_{s},
\int_{0}^{s}k\left(t,\tau,y_{\tau}+\widetilde{\phi}_{\tau}\right)d\tau\right)ds\right\|_{X}\\
&+&
\left\|\lim\limits_{\lambda\longrightarrow+\infty}\int_{t_{1}-\epsilon}^{t_{1}}\left(S'(t_{2}-s)-S'(t_{1}-s)\right)\right.\\
&\times&\left.B_{\lambda}f\left(s,y_{s}+\widetilde{\phi}_{s},
\int_{0}^{s}k\left(t,\tau,y_{\tau}+\widetilde{\phi}_{\tau}\right)d\tau\right)ds\right\|_{X}
\\
&+&
\left\|\lim\limits_{\lambda\longrightarrow+\infty}\int_{t_{1}}^{t_{2}}S'(t_{2}-s)B_{\lambda}f\left(s,y_{s}+\widetilde{\phi}_{s},
\int_{0}^{s}k\left(t,\tau,y_{\tau}+\widetilde{\phi}_{\tau}\right)d\tau\right)ds\right\|_{X}
\\
&+&\left\|\sum\limits_{0<t_{k}<t_{1}}\left(S'(t_{1}-t_{k})-S'(t_{2}-t_{k})\right)I_{k}\left(y(t_{k}^{-})+\widetilde{\phi}(t_{k}^{-})\right)\right\|_{X}\\
&+&\left\|\sum\limits_{t_{1}<t_{k}<t_{2}}S'(t_{2}-t_{k})I_{k}\left(y(t_{k}^{-})+\widetilde{\phi}(t_{k}^{-})\right)\right\|_{X}\\
&\leq&
\lim\limits_{\lambda\longrightarrow+\infty}\left(\int_{0}^{t_{1}-\epsilon}\left\|\left(S'(t_{2}-s)-S'(t_{1}-s)\right)\right\|_{X}\left\|B_{\lambda}\right\|_{X}\right.\\
&\times& \alpha_{q'\overline{\omega}}(s)ds
+\int_{t_{1}-\epsilon}^{t_{1}}\left\|\left(S'(t_{2}-s)-S'(t_{1}-s)\right)\right\|_{X}\left\|B_{\lambda}\right\|_{X}\alpha_{q'\overline{\omega}}(s)ds\\
&+&\left.\int_{t_{1}}^{t_{2}}\left\|S'(t_{2}-s)\right\|_{X}\left\|B_{\lambda}\right\|_{X}\alpha_{q'\overline{\omega}}(s)ds\right)\\
&+&\sum\limits_{0<t_{k}<t_{1}}\left(S'(t_{1}-t_{k})-S'(t_{2}-t_{k})\right)L_{k}(l^{-1}q')\\
&+& \sum\limits_{t_{1}<t_{k}<t_{2}}S'(t_{2}-t_{k})L_{k}(l^{-1}q').
\end{eqnarray*}
Since $S'(t)$ is compact for $t>0$ and hence continuous in the
uniform operator topology, $\alpha_{q'\overline{\omega}}\in L^{1}$,
and the right hand side of is independent of $y\in B_{q}$, we obtain
$\left\|\psi_{2}(y)(t_{1})-\psi_{2}(y)(t_{2})\right\|\longrightarrow0$
as $\left(t_{1}-t_{2}\right)\longrightarrow0$ with $\epsilon>0$
sufficiently small. This prove that $\psi_{2}$ maps $B_{q}$ into an
equicontinuous family of functions.\\

\begin{lemma} $\psi_{2}$ maps $B_{q}$ into a precompact set in
$X$.
\end{lemma}
Proof.\\ Let $0<t\leq b$ be fixed and $\epsilon$ be
a real number satisfying $0<\epsilon<t$. For $y\in B_{q}$, we define
\begin{eqnarray*}
\left(\psi_{2}y\right)(t)&=&\lim\limits_{\lambda\longrightarrow+\infty}\int_{0}^{t-\epsilon}S'(t-s)B_{\lambda}f\left(s,y_{s}+\widetilde{\phi}_{s},\int_{0}^{s}
k\left(t,\tau,y_{\tau}+\widetilde{\phi}_{\tau}\right)d\tau\right)ds\\
&+&\sum\limits_{0<t_{k}< t}
S'(t-t_{k})\delta(t-t_{k})I_{k}\left(y(t_{k}^{-})+\widetilde{\phi}(t_{k}^{-})\right)\\
&=&S'(\epsilon)\lim\limits_{\lambda\longrightarrow+\infty}\int_{0}^{t-\epsilon}S'(t-\epsilon-s)\\
&\times& B_{\lambda}f\left(s,y_{s}+\widetilde{\phi}_{s},
\int_{0}^{s}k\left(t,\tau,y_{\tau}+\widetilde{\phi}_{\tau}\right)d\tau\right)ds\\
&+&\sum\limits_{0<t_{k}< t}
S'(t-t_{k})\delta(t-t_{k})I_{k}\left(y(t_{k}^{-})+\widetilde{\phi}(t_{k}^{-})\right).
\end{eqnarray*}
Note that,
\begin{equation}\label{estimate}
\left\|\lim\limits_{\lambda\longrightarrow+\infty}\int_{0}^{t-\epsilon}S'(t-\epsilon-s)B_{\lambda}f\left(s,y_{s}+\widetilde{\phi}_{s},
\int_{0}^{s}k(s,\tau,y_{\tau}+\widetilde{\phi}_{\tau})d\tau\right)ds\right\|_{X}\leq
M\overline{M}\int_{0}^{t-\epsilon}\alpha_{q'\overline{\omega}}(s)ds,
\end{equation}
Using the estimations (\ref{psi2}), (\ref{estimate}) and by the
compactness of $S'(t)$ ($t>0$), we obtain that the set
$Z_{\epsilon}(t)=\left\{\left(\psi_{2\epsilon}y\right)(t):y\in
B_{q}\right\}$ is relative compact in $X$ for every $\epsilon$,
$0<\epsilon<t$. Thus, for every $q\in B_{q}$, we have
\begin{eqnarray*}
\left\|\left(\psi_{2}y\right)(t)-\left(\psi_{2\epsilon}y\right)(t)\right\|_{X}&\leq&\left\|\lim\limits_{\lambda\longrightarrow+\infty}
\int_{t-\epsilon}^{t}S'(t-s)B_{\lambda}f\left(s,y_{s}+\widetilde{\phi}_{s},\int_{0}^{t}k\left(t,\tau,y_{\tau}+\widetilde{\phi}_{\tau}\right)d\tau\right)ds\right\|_{X}\\
&\leq&\int_{t-\epsilon}^{t}M\overline{M}\alpha_{q'\overline{\omega}}(s)ds.
\end{eqnarray*}
Therefore
$\left\|\left(\psi_{2}y\right)(t)-\left(\psi_{2\epsilon}y\right)(t)\right\|\longrightarrow0$
as $\epsilon\longrightarrow0^{+}$, and hence there are precompact
sets arbitrarily close to the set
$\left\{\left(\psi_{2\epsilon}y\right)(t):y\in B_{q}\right\}$, thus
this set is also relative compact in $X$.\\

\begin{lemma}\label{7}
$\psi_{2}:\mathrm{BM}_{h}''\longrightarrow\mathrm{BM}_{h}''$ is
continuous.
\end{lemma}
Proof.\\ Let
$\left\{y^{(n)}\right\}_{n=0}^{+\infty}\subseteq\mathrm{BM}_{h}''$,
with $y^{(n)}\longrightarrow y$ in $\mathrm{BM}_{h}''$. Therefore
there exists $q>0$ such that $\left\|y^{(n)}(t)\right\|_{X}\leq q$
for
all $n\in \mathbb{N}$ and $t\in J$, then $y^{(n)},y\in B_{q}$.\\
Using Lemma \ref{encadre},
$\left\|y_{t}^{(n)}+\widetilde{\phi}_{t}\right\|_{\mathrm{BM}_{h}}\leq
q'$, for all $t\in J$, and by using the hypothesis ($H_{6}$), we
obtain \\

(i) For each $k$, ($k=1,...,m$), $I_{k}$ is continuous,\\

(ii)
$f\left(t,y_{t}^{(n)}+\widetilde{\phi}_{t},\int_{0}^{t}k\left(t,s,y_{s}^{(n)}+\widetilde{\phi}_{s}\right)ds\right)\longrightarrow
f\left(t,y_{t}+\widetilde{\phi}_{t},\int_{0}^{t}k\left(t,s,y_{s}+\widetilde{\phi}_{s}\right)ds\right)$.
By the inequality
$$
\left\|f\left(t,y_{t}^{(n)}+\widetilde{\phi}_{t},\int_{0}^{t}k\left(t,s,y_{s}^{(n)}+\widetilde{\phi}_{s}\right)ds\right)-
f\left(t,y_{t}+\widetilde{\phi}_{t},\int_{0}^{t}k\left(t,s,y_{s}+\widetilde{\phi}_{s}\right)ds\right)\right\|_{X}\leq2\alpha_{q'\overline{\omega}}(t),
$$
and using the dominated convergence theorem,
\begin{eqnarray*}
\left\|\left(\psi_{2\epsilon}y^{(n)}\right)(t)-\left(\psi_{2\epsilon}y\right)(t)\right\|_{X}&=&
\left\|\lim\limits_{\lambda\longrightarrow+\infty}\int_{0}^{t}S'(t-s)B_{\lambda}\left[
f\left(s,y_{s}^{(n)}+\widetilde{\phi}_{s},\int_{0}^{s}k\left(t,\tau,y_{\tau}^{(n)}+\widetilde{\phi}_{\tau}\right)d\tau\right)\right.\right.\\
&-&\left.\left.f\left(
s,y_{s}+\widetilde{\phi}_{s},\int_{0}^{s}k\left(t,\tau,y_{\tau}+\widetilde{\phi}_{\tau}\right)\right)\right]\right\|_{X}\\
&+& \left\|\sum\limits_{0<t_{k}<t}S'(t-t_{k})\left[
\delta(t-t_{k})I_{k}(y^{(n)}(t_{k}^{-})+\widetilde{\phi}(t_{k}^{-}))\right]\right\|_{X}\\
&\leq&
M\overline{M}\int_{0}^{t}\left\|f\left(t,y_{t}^{(n)}+\widetilde{\phi}_{t},\int_{0}^{s}k\left(t,\sigma,y_{\sigma}^{(n)}+\widetilde{\phi}_{\sigma}\right)d\sigma\right)\right.\\
&-&
\left.f\left(t,y_{t}+\widetilde{\phi}_{t},\int_{0}^{s}k\left(\sigma,y_{\sigma}+\widetilde{\phi}_{\sigma}\right)d\sigma\right)\right\|_{X} ds\\
&+&
M\sum\limits_{0<t_{k}<t}\left\|I_{k}\left(y^{(n)}(t_{k}^{-})+\widetilde{\phi}(t_{k}^{-})\right)-I_{k}\left(y(t_{k}^{-})+\widetilde{\phi}(t_{k}^{-})\right)\right\|_{X},
\end{eqnarray*}
thus
$$\left\|\left(\psi_{2\epsilon}y^{(n)}\right)(t)-\left(\psi_{2\epsilon}y\right)(t)\right\|_{X}\longrightarrow0 \,\ \text{as} \,\ n\longrightarrow+\infty.$$
Therefore,
$\left\|\psi_{2}y^{(n)}-\psi_{2}y\right\|_{b}=\sup\limits_{t\in
J}\left\|\psi_{2}y^{(n)}(t)-\psi_{2}y(t)\right\|_{X}\longrightarrow0$,
and hence $\psi_{2}$ is continuous.\\

\begin{theorem}\label{th1} Assume that assumptions ($H_{1}$)-($H_{6}$) are fulfilled,
$\phi\in\mathrm{BM}_{h}$ and $\phi(0)-g(0,\phi,0)\in
\overline{D(A)}$. Then (\ref{elsevier1})-(\ref{elsevier111}) has at
least one piecewise pseudo almost periodic solution on $(-\infty,b]$
provided that
$$
K:=l(G_{1}+bG_{2}H)<1,
$$
$$
(G_{1}+G_{2})\overline{\omega}l+M\overline{M}\overline{\omega}l\mu+M\sum\limits_{k=1}^{m}\frac{L_{k}(\sigma)}{\sigma}<1,\,\
\text{where}\,\ \overline{\omega}=\max\{1,p(t),q(t)\}.
$$
\end{theorem}
Proof.\\ By using Sadovskii's fixed point theorem we prove that the
operator $\Gamma$ defined by
\begin{equation*}\Gamma(x)(t)=\left\{
\begin{array}{ll}
\phi(t), \,\ \text{if}\,\ t\in(-\infty,0],\\
S'(t)\left[\phi(0)-g(0,\phi,0)\right]+g\left(t,x_{t},\int_{0}^{t}h(t,s,x_{s})ds\right)\\
+\lim\limits_{\lambda\longrightarrow+\infty}\int_{0}^{t}S'\left(t-s\right)\mathrm{B}_{\lambda}f(s,x_{s},\int_{0}^{s}k(t,\tau,x_{\tau})d\tau)ds\\
+\sum\limits_{0<t_{k}<t}S'\left(t-t_{k}\right)\delta(t-t_{k})I_{k}\left(x(t_{k}^{-})\right),
\,\ t\in J,
\end{array}\right.
\end{equation*}has fixed point $x(.)$. This fixed point is then a PPAP solution of the system (\ref{elsevier1})-(\ref{elsevier111}).\\

For $\phi\in\mathrm{BM}_{h}$, define $\widetilde{\phi}$ by
\begin{equation*}
\widetilde{\phi}(t)=\left\{
\begin{array}{ll}
\phi(t),\,\ t\in(-\infty,0],\\
S'(t)\phi(0), \,\ t\in J=[0,b],
\end{array}
\right.
\end{equation*}
then $\widetilde{\phi}\in\mathrm{BM}_{h}'$. Let
$x(t)=y(t)+\widetilde{\phi}(t)$, $-\infty<t\leq b$. We remark that
the shape of $x$ satisfies (\ref{elsevier2}) if and only if $y$
satisfies $y_{0}=0$ and
\begin{eqnarray*}
y(t)&=&-S'(t)g(0,\phi,0)+g\left(t,y_{t}+\widetilde{\phi}(t),\int_{0}^{t}h\left(t,s,y_{s}+\widetilde{\phi}_{s}\right)ds\right)\\
&+&
\lim\limits_{\lambda\longrightarrow+\infty}\int_{0}^{t}S'(t-s)B_{\lambda}f\left(s,y_{s}+\widetilde{\phi}_{s},\int_{0}^{s}k\left(t,\tau,y_{\tau}+\widetilde{\phi}_{\tau}\right)d\tau\right)ds\\
&+& \sum\limits_{0<t<
t_{k}}S'(t-t_{k})\delta(t-t_{k})I_{k}\left(y(t_{k}^{-})+\widetilde{\phi}(t_{k}^{-})\right),
\,\ t\in J.
\end{eqnarray*}

From Lemma \ref{1}-Lemma \ref{7} and the Arzela-Ascoli theorem it
follows that $\psi_{2}:B_{q}\longrightarrow B_{q}$ is completely
continuous, thus, we have proved that $\psi$ is condensing operator
on $B_{q}$. Using the Sadovskii fixed point theorem, $\psi$ has a
fixed point $\theta(.)$ in $B_{q}$. Let
$x(t)=\theta(t)+\widetilde{\phi}(t)$, $t\in(-r,b]$, for all $r>0$.
Then $x$ is an PPAP solution of
(\ref{elsevier1})-(\ref{elsevier111}). This achieve the proof of
Theorem \ref{th1}.

\section{Application}

Let us consider the following partial neutral mixed integro
differential equation with infinite delays, presented by the system:

\begin{eqnarray}\label{elsevier12}
&&\frac{d}{dt}\left[z(t,x)-G\left(t,\int_{-\infty}^{t}P_{1}(s-t)z(s,x)ds,
\int_{0}^{t}\int_{-\infty}^{s}P_{2}(s,x,\tau-s)Q_{1}(z(\tau,x))d\tau
ds\right)\right]\nonumber\\&=&
 \frac{\partial^{2}}{\partial x^{2}}\left[z(t,x)-G\left(t,\int_{-\infty}^{t}P_{1}(s-t)z(s,x)ds,
\int_{0}^{t}\int_{-\infty}^{s}P_{2}(s,x,\tau-s)Q_{1}(z(\tau,x))d\tau
ds\right)\right]\nonumber\\
&+&
F\left(t,\int_{-\infty}^{0}P_{3}(s-t)z(s,x)ds,\int_{0}^{t}\int_{-\infty}^{s}P_{4}(s,z,\tau-s)Q_{2}(z(\tau,x))d\tau
ds \right),\\
&& \text{if}\,\ (t,x)\in [0,b]\backslash
t_{k}\times[0,\pi],
\nonumber\\
 && \label{elsevier122}
z(t,0)=z(t,\pi)=0, \,\ 0\leq t\leq b,\\
&& z(t_{k}^{+},x)-z(t_{k}^{-},x)=I_{k}(z(t_{k}^{-},x)), \,\ \text{if} \,\ t=t_{k}\,\ \text{and}\,\ x\in[0,\pi],\\
 &&\label{elsevier1222} z(t,x)=\phi(t,x),\,\ t\in [-r,0], \,\
x\in[0,\pi], \,\ r>0,\,\ k=1,...,m, \nonumber
\end{eqnarray}

where
$F:[0,b]\times\mathbb{R}\times\mathbb{R}\longrightarrow\mathbb{R}$
and the functions
$P_{2},P_{4}:\Delta\times\mathbb{R}\longrightarrow\mathbb{R}$ are
continuous
and strongly mesurable functions.\\

Let $X=L^{2}[0,\pi]$ with the norm $\|\cdot\|_{L^{2}}$, define the
operator
$$
A:D(A)\subset X\longrightarrow X \,\ \text{by}\,\
Aw=w'',
$$
with domain
$$
D(A)=\{w\in X: w,w'\,\ \text{are absolutely continuous}\,\ w''\in X,
\,\ w(0)=w(\pi)=0\},
$$
then we have
$$
\overline{D(A)}=\{w\in X: w,w'\,\ \text{are absolutely
continuous}\,\ w(0)=w(\pi)=0\}\neq X.
$$
It is well known that (see \cite{prato}) $A$ satisfies the following
properties:\\

(i) $\rho(A)\supseteq (0,+\infty)$.\\

(ii) $\left\|(\lambda I-A)^{-1}\right\|\leq \frac{1}{\lambda}$, for
all $\lambda>0$.\\

Then $Aw=\sum\limits_{n=1}^{+\infty}n^{2}<w,w_{n}> w_{n}$, and $w\in
D(A)$, where $w_{n}(x)=\sqrt{\frac{2}{\Pi}}\sin(nx)$, $n=1,2,...$ is
the orthogonal set of eigen vectors of $A$. It is well known that
$A$ is the infinitesimal generator of an analytic semigroup $S(t)$,
$t\geq0$, and is given by
$S(t)w=\sum\limits_{n=1}^{+\infty}e^{-n^{2}t}<w,w_{n}>w_{n}$, $w\in
X$.

Since the analytic semigroup $\{S'(t)\}_{t\geq0}$ is compact
\cite{pazy}, there exists constants $M>0$ such that $\|S'(t)\|\leq
M$.

Further for every $w\in X$,
$(-A)^{-\frac{1}{2}}w=\sum\limits_{n=1}^{+\infty}\frac{1}{n}<w,w_{n}>w_{n}$,
and
$$\|(-A)\|^{-\frac{1}{2}}=1.$$ The operator $(-A)^{\frac{1}{2}}$
is given by
$$
(-A)^{\frac{1}{2}}=\sum\limits_{n=1}^{+\infty}n<w,w_{n}>w_{n},
$$
with the domain $D((-A)^{\frac{1}{2}})=\left\{w\in X: \,\
\sum\limits_{n=1}^{+\infty}n<w,w_{n}>w_{n}\in X\right\}$. It follows
that $\|S(t)\|\leq1$. Let $h(s)=\exp(2s)$, $s<0$; then
$l=\int_{-\infty}^{0}h(s)ds=\frac{1}{2}$, and define
$$
\|\phi\|_{\mathrm{BM}_{h}}=\int_{-\infty}^{0}h(s)\|\phi(\theta)\|_{L^{2}}ds.
$$

Hence for $(t,\phi)\in[0,b]\times \mathrm{BM}_{h}$, where
$\phi(\theta)(x)=\phi(\theta,x)$,
$(\theta,x)\in(-\infty,0]\times[0,\pi]$. Set
$$z(t)(y)=z(t,y);$$
$$
k(t,s,\psi)(y)=\int_{-\infty}^{s}P_{4}(t,y,\theta)Q_{2}(\psi(\theta)y)d\theta;
$$
$$
h(t,s,\psi)(y)=\int_{-\infty}^{s}P_{2}(t,y,\theta)Q_{1}(\psi(\theta)y)d\theta;
$$
$$
g\left(t,\psi,\int_{0}^{t}h(t,s,\psi)ds\right)(y)=G\left(t,\int_{-\infty}^{0}P_{1}(\theta)\psi(\theta)(y)d\theta,\int_{0}^{t}h(t,s,\psi)(y)ds\right);
$$
$$
f\left(t,\psi,\int_{0}^{t}k(t,s,\psi)ds\right)(y)=F\left(t,\int_{-\infty}^{0}P_{3}(\theta)\phi(\theta)(y)d\theta,\int_{0}^{t}k(t,s,\psi)(y)
ds \right).
$$

The above partial differential system
(\ref{elsevier12})-(\ref{elsevier1222}) can be formulated as an
abstract form as the system (\ref{elsevier1})-(\ref{elsevier111}).
Suppose further that:\\

($a_{1}$) The function $P_{2}(t,y,\theta)$ is continuous in
$[0,b]\times[0,\pi]\times(-\infty,0]$ and $P_{2}(t,y,\theta)\geq0$, $\int_{0}^{t}\int_{-\infty}^{s}P_{2}(t,y,\theta)d\theta ds=\tilde{p}_{2}(t,y)<\infty$.\\

($a_{2}$) The function $P_{4}(t,y,\theta)$ is continuous in
$[0,b]\times[0,\pi]\times(-\infty,0]$ and $P_{4}(t,y,\theta)\geq0$, $\int_{0}^{t}\int_{-\infty}^{s}P_{4}(t,y,\theta)d\theta ds=\tilde{p}_{4}(t,y)<\infty$.\\

($a_{3}$) The function $Q_{i}(\cdot)$ is continuous and for
$(\theta,y)\in(-\infty,0]\times[0,\pi]$,
$$
0\leq Q_{i}(\psi(\theta)y)\leq \vartheta_{i}
\left(\int_{-\infty}^{0}e^{2s}\|\psi(s,\cdot)\|_{L^{2}}ds\right),\,\
i=1,2,
$$
where $\vartheta_{i}: [0,+\infty)\longrightarrow (0,+\infty)$ is a
continuous and nondecreasing function.\\

Now we can see that
\begin{eqnarray*}
\left\|\int_{0}^{t}h(t,s,\phi)(y)ds\right\|_{L^{2}}&=&
\left[\int_{0}^{\pi}\left(\int_{0}^{t}\int_{-\infty}^{s}P_{2}(t,y,\theta)Q(\phi(\theta)(y))d\theta ds\right)^{2}dy\right]^{\frac{1}{2}}\\
&\leq&
\left[\int_{0}^{\pi}\left(\int_{0}^{t}\int_{-\infty}^{s}P_{2}(t,y,\theta)\vartheta_{1}\left(\int_{-\infty}^{0}e^{2s}\|\phi(s)(\cdot)\|_{L^{2}}ds\right)d\theta ds\right)^{2}dy\right]^{\frac{1}{2}}\\
&\leq&\left[\int_{0}^{\pi}\left(\int_{0}^{t}\int_{-\infty}^{s}P_{2}(t,y,\theta)\vartheta_{1}\left(\int_{-\infty}^{0}e^{2s}\sup\limits_{s\in[\theta,0]}\|\phi(s)\|_{L^{2}}ds\right)d\theta ds\right)^{2}dy\right]^{\frac{1}{2}}\\
&\leq&\left[\int_{0}^{\pi}\left(\int_{0}^{t}\int_{-\infty}^{s}P_{2}(t,y,\theta)d\theta ds\right)^{2}dy\right]^{\frac{1}{2}}\vartheta_{1}(\|\phi\|_{\mathrm{BM}_{h}})\\
&\leq&\left[\int_{0}^{\pi}(\tilde{p}_{2}(t,y))^{2}dy\right]^{\frac{1}{2}}\vartheta_{1}(\|\phi\|_{\mathrm{BM}_{h}})\triangleq p_{2}(t) \vartheta_{1}(\|\phi\|_{\mathrm{BM}_{h}}),\\
\end{eqnarray*}
and
\begin{eqnarray*}
\|k(t,s,\phi)(y)\|_{L^{2}}&=&
\left[\int_{0}^{\pi}\left(\int_{0}^{t}\int_{-\infty}^{s}P_{4}(t,y,\theta)Q(\phi(\theta)(y))d\theta ds\right)^{2}dy\right]^{\frac{1}{2}}\\
&\leq&
\left[\int_{0}^{\pi}\left(\int_{0}^{t}\int_{-\infty}^{s}P_{4}(t,y,\theta)\vartheta_{2}\left(\int_{-\infty}^{0}e^{2s}\|\phi(s)(\cdot)\|_{L^{2}}ds\right)d\theta
ds\right)^{2}dy\right]^{\frac{1}{2}}
\\
&\leq&\left[\int_{0}^{\pi}\left(\int_{0}^{t}\int_{-\infty}^{s}P_{4}(t,y,\theta)\vartheta_{2}\left(\int_{-\infty}^{0}e^{2s}\sup\limits_{s\in[\theta,0]}\|\phi(s)\|_{L^{2}}ds\right)d\theta ds\right)^{2}dy\right]^{\frac{1}{2}}\\
&\leq&\left[\int_{0}^{\pi}\left(\int_{0}^{t}\int_{-\infty}^{s}P_{4}(t,y,\theta)d\theta ds\right)^{2}dy\right]^{\frac{1}{2}}\vartheta_{2}(\|\phi\|_{\mathrm{BM}_{h}})\\
&\leq&\left[\int_{0}^{\pi}(\tilde{p}_{4}(t,y))^{2}dy\right]^{\frac{1}{2}}\vartheta_{2}(\|\phi\|_{\mathrm{BM}_{h}})\triangleq
p_{4}(t) \vartheta_{2}(\|\phi\|_{\mathrm{BM}_{h}}).
\end{eqnarray*}

By imposing suitable conditions on the above defined pseudo almost
periodic functions $F$, $G$ and $I_{k}$, to verify the assumptions
of Theorem \ref{th1}, we conclude that the system
(\ref{elsevier1})-(\ref{elsevier111}) has at least one PPAP
solution.

\section{Conclusion}
We have considered an impulsive integro-differential system with
infinite delays. This research serves as a first step for the
application of the Sadovskii's theorem in investigation of the
existence of PPAP solutions for such systems. After introducing an
abstract phase space $\mathrm{BM}_{h}$ (in Section 2). Some result
are improved and generalized. The model investigated in this work
can be regarded as a natural continuation of the models studied in
(\cite{abada}, \cite{acta1}, \cite{aacta},
\cite{neuro}, \cite{chang}, \cite{han}, \cite{kucche}, \cite{shan},
\cite{sstamov}, \cite{zabut}). An example is also included to
illustrate the importance of the results obtained. The demonstrated
techniques can be applied in studying of qualitative properties of
many practical problems in diverse domain. Besides, this techniques
can been extended to investigate the controllability (see
\cite{sak}) of the first order impulsive functional
integro-differential with infinite delay. In the best of our knowledge, we can applied this method for the model investigated in (\cite{supp1}, \cite{supp2}, \cite{supp3}, \cite{supp4}).

\end{document}